# SPECTRAL THEORY FOR STRUCTURED PERTURBATIONS OF LINEAR OPERATORS

MARTIN ADLER AND KLAUS-JOCHEN ENGEL

ABSTRACT. We characterize the spectrum (and its parts) of operators which can be represented as "$G = A + BC$" for a "simpler" operator $A$ and a structured perturbation $BC$. The interest in this kind of perturbations is motivated, e.g., by perturbations of the domain of an operator $A$ but also arises in the theory of closed-loop systems in control theory. In many cases our results yield the spectral values of $G$ as zeros of a *"characteristic equation"*.

## 1. INTRODUCTION

The *spectrum*
$$\sigma(G) := \big\{\lambda \in \mathbb{C} : \lambda - G \text{ is not invertible in } \mathcal{L}(X)\big\}$$
as a subset of $\mathbb{C}$, and its finer subdivisions (cf. Definition A.6) reflect much information about a (possibly unbounded) linear operator $G : D(G) \subset X \to X$ on a (generally infinite dimensional) Banach space $X$. Here we only mention (for details see [7, Chap. V]) that for generators $G$ of strongly continuous semigroups the location of $\sigma(G)$ in the complex plane determines to a great part the asymptotic behavior of the solutions of the associated abstract Cauchy problem

(ACP) $\qquad \begin{cases} \frac{d}{dt} x(t) = Gx(t), & t \geq 0, \\ x(0) = x_0. \end{cases}$

However, in many applications it is difficult to determine $\sigma(G)$ by direct computations. One approach to overcome this difficulty is to split $G$ into a suitable sum "$G = A + P$" of a well understood operator $A$ and a perturbation $P$. Then one tries to characterize spectral values of $G$ by "simple" conditions involving the operators $A$ and $P$.

In this paper we elaborate this idea for "structured perturbations", i.e., perturbations which can be written as a product $P = BC$. First we setup our general framework. For a summary of our notation we refer to Appendix A.1.

*Assumptions* 1.1. We assume that
 (i) $U$, $X$, $Z$ and $Z_{-1}$ are Banach spaces such that $Z \hookrightarrow X \hookrightarrow Z_{-1}$;
 (ii) $A^Z : Z \subset Z_{-1} \to Z_{-1}$ is a linear operator satisfying $\rho(A^Z) \neq \emptyset$ and $A^Z \in \mathcal{L}(Z, Z_{-1})$;
 (iii) $B \in \mathcal{L}(U, Z_{-1})$ and $C \in \mathcal{L}(Z, U)$.

We are then interested in the operator $G = A_{BC} : D(A_{BC}) \subseteq X \to X$ given by

(1.1) $\qquad A_{BC} := (A^Z + BC)|_X \quad \text{with domain} \quad D(A_{BC}) := \big\{x \in Z : (A^Z + BC)x \in X\big\},$

where the sum is initially taken in $Z_{-1}$. This setting is summarized in Diagram 1. For the main cases fitting into this setup see Sections 3–4 and Remark A.5.(ii).

If we define the operator $A : D(A) \subseteq X \to X$ by

(1.2) $\qquad A := A^Z|_X \quad \text{with domain} \quad D(A) := \big\{x \in Z : A^Z x \in X\big\},$

then we can consider $A_{BC}$ as a perturbation of $A$, where Assumption 1.1.(iii) limits the unboundedness of the structured perturbation $P := BC \in \mathcal{L}(Z, Z_{-1})$. Now Lemma A.7.(vi)–(vii) gives the following result.

---

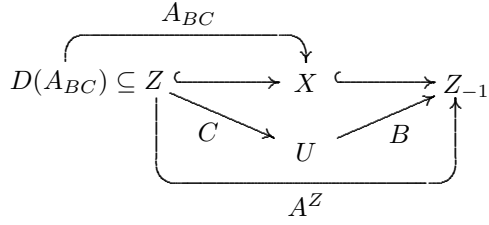

Diagram 1: The operators defining $A_{BC}$ in (1.1).

**Lemma 1.2.** *Let $A$ be given by (1.2). Then*

$$(1.3) \qquad \rho(A) = \rho(A^Z) \qquad \text{and} \qquad R(\lambda, A) = R(\lambda, A^Z)|_X \quad \text{for all } \lambda \in \rho(A).$$

In [1, 2] we studied in detail the generator property of $A_{BC}$, cf. Remark A.5.(ii). In the present paper we characterize the various types of spectral values in terms of an operator $\Delta_W(\lambda)$ defined on some, in general "smaller", space $W$. In particular, if $\dim(W) < +\infty$, this yields the spectral values in $\sigma(A_{BC}) \cap \rho(A)$ as zeros of a (nonlinear) *characteristic equation*, cf. (2.2) in Theorem 2.3. Our interest in operators $A_{BC}$ given by (1.1) is, among other, motivated by perturbations of the domain of operators in the spirit of Greiner, cf. [8]. For this reason, in Subsection 3 we first apply our main abstract result Theorem 2.3 to this generic situation. The usefulness of our approach is then demonstrated by a series of further concrete examples in Section 4.

In Appendix A we summarize the notation, give a short introduction to the extrapolation of spaces and operators, prove some results concerning the spectral theory of parts of operators and present so-called Schur complements for operator matrices needed for our approach.

We mention that related problems have already been studied by, e.g., Salamon, Weiss–Xu and Curtain–Jacob in the context of closed loop systems in control theory, cf. [16, Lem. 4.4], [19, Thms. 1.1 & 1.2] and [6, Thm. 6.2]. Our Theorem 2.3 generalizes these results since it does not rely on admissibility conditions for the operators $B$, $C$, a Hilbert space structure of $X$ or the closedness or a dense domain of $A_{BC}$. Moreover, we study not only the spectrum and point spectrum but characterize also other parts as the approximate point-, continuous-, residual- and essential spectrum of $A_{BC}$.

## 2. Spectral Theory for $A_{BC}$

In this section we investigate the spectrum of the perturbed operator $A_{BC}$ from (1.1) limiting the hypotheses to Assumption 1.1. In particular, we do *not* assume $X$ to be a Hilbert space or $A_{BC}$ to be closed or densely defined. Moreover, we do not impose any kind of "admissibility condition" on the triple $(A, B, C)$. The proof of our main result Theorem 2.3 is based on two ingredients: spectral properties of the part of an operator in a subspace, see Appendix A.3, and Schur complements for operator matrices, cf. Appendix A.4.

To start our investigations, we define the operator $A_{BC}^Z : D(A_{BC}^Z) \subseteq Z_{-1} \to Z_{-1}$ by

$$A_{BC}^Z := A^Z + BC \quad \text{with domain} \quad D(A_{BC}^Z) := Z,$$

for which the following holds.

**Lemma 2.1.** *We have $A_{BC}^Z \in \mathcal{L}(Z, Z_{-1})$. Moreover, if $A_{BC}^Z$ is closed, e.g., $\rho(A_{BC}^Z) \neq \emptyset$, then the norm of $Z$ and the graph norm of $A_{BC}^Z$ are equivalent on $Z$, i.e.,*

$$(2.1) \qquad \|\bullet\|_Z \simeq \|\bullet\|_{A_{BC}^Z},$$

*where $\|x\|_{A_{BC}^Z} := \|x\|_{Z_{-1}} + \|A_{BC}^Z x\|_{Z_{-1}}$ for $x \in Z$. In other words, $Z \simeq (Z_{-1})_1^{A_{BC}^Z}$.*

*Proof.* By assumption $BC \in \mathcal{L}(Z, Z_{-1})$ and $A^Z \in \mathcal{L}(Z, Z_{-1})$, hence $A_{BC}^Z \in \mathcal{L}(Z, Z_{-1})$. Moreover, $Z \hookrightarrow Z_{-1}$ and therefore $\|\bullet\|_Z$ is finer than $\|\bullet\|_{A_{BC}^Z}$. Now, if $A_{BC}^Z$ is closed then $(Z, \|\bullet\|_{A_{BC}^Z})$ is a Banach space and the equivalence in (2.1) follows from the open mapping theorem. □

The following operators will be the main tool in the sequel. Recall that $\rho(A) = \rho(A^Z)$ by (1.3).

**Definition 2.2.** For $\lambda \in \rho(A)$ define the operators
$$\Delta_U(\lambda) := CR(\lambda, A^Z)B \in \mathcal{L}(U) \quad \text{and} \quad \Delta_Z(\lambda) := R(\lambda, A^Z)BC \in \mathcal{L}(Z).$$

Here the boundedness of $\Delta_W(\lambda)$ for $W \in \{U, Z\}$ follows from Assumption 1.1, the closed graph theorem and the resolvent equation. Using these operators, the spectral values of $A_{BC}$ can be characterized in the following way. For the notions concerning the finer division of the spectrum, see Definition A.6.

**Theorem 2.3.** *Let $\lambda \in \rho(A)$ and $W \in \{U, Z\}$.*

(a) *The following spectral inclusions always hold.*
$$\sigma(A_{BC}) \subseteq \sigma(A_{BC}^Z), \qquad \sigma_p(A_{BC}) = \sigma_p(A_{BC}^Z),$$
$$\sigma_a(A_{BC}) \subseteq \sigma_a(A_{BC}^Z), \qquad \sigma_{ess}(A_{BC}) \subseteq \sigma_{ess}(A_{BC}^Z).$$

*If $D(A) + \mathrm{rg}(Id_Z - \Delta_Z(\nu))$ is dense in $Z$ for some (hence all) $\nu \in \rho(A)$ then also*
$$\sigma_c(A_{BC}) \subseteq \sigma_c(A_{BC}^Z), \qquad \sigma_r(A_{BC}) \supseteq \sigma_r(A_{BC}^Z).$$

(b) *The following spectral characterizations always hold.*
$$\lambda \in \sigma(A_{BC}^Z) \quad \Longleftrightarrow \quad 1 \in \sigma(\Delta_W(\lambda)),$$
$$\lambda \in \sigma_*(A_{BC}^Z) \quad \Longleftrightarrow \quad 1 \in \sigma_*(\Delta_W(\lambda))$$

*for all $* \in \{p, a, r, c, ess\}$. Moreover,*
$$\lambda \in \sigma(A_{BC}) \quad \Longrightarrow \quad 1 \in \sigma(\Delta_W(\lambda)),$$
$$\lambda \in \sigma_p(A_{BC}) \quad \Longleftrightarrow \quad 1 \in \sigma_p(\Delta_W(\lambda)),$$
$$\lambda \in \sigma_a(A_{BC}) \quad \Longrightarrow \quad 1 \in \sigma_a(\Delta_W(\lambda)),$$
$$\lambda \in \sigma_{ess}(A_{BC}) \quad \Longrightarrow \quad 1 \in \sigma_{ess}(\Delta_W(\lambda)).$$

*If $D(A) + \mathrm{rg}(Id_Z - \Delta_Z(\nu))$ is dense in $Z$ for some (hence all) $\nu \in \rho(A)$ then also*
$$\lambda \in \sigma_c(A_{BC}) \quad \Longrightarrow \quad 1 \in \sigma_c(\Delta_W(\lambda)),$$
$$\lambda \in \sigma_r(A_{BC}) \quad \Longleftarrow \quad 1 \in \sigma_r(\Delta_W(\lambda)).$$

(c) *If $\Delta_W(\lambda) \in \mathcal{L}(W)$ is compact, then*
$$\lambda \in \sigma(A_{BC}) \quad \Longleftrightarrow \quad \lambda \in \sigma_p(A_{BC}) \quad \Longleftrightarrow \quad 1 \in \sigma_p(\Delta_W(\lambda)).$$

*In particular, if $\dim(U) < \infty$, then*

(2.2) $\qquad \lambda \in \sigma(A_{BC}) \quad \Longleftrightarrow \quad \lambda \in \sigma_p(A_{BC}) \quad \Longleftrightarrow \quad \det(Id_U - \Delta_U(\lambda)) = 0.$

(d) *If the condition*

(2.3) $\qquad 1 \in \rho(\Delta_W(\nu)) \text{ for some } \nu \in \rho(A) \qquad \text{or, equivalently,} \qquad \rho(A) \cap \rho(A_{BC}^Z) \neq \emptyset$

*holds, then for all $* \in \{a, r, c, ess\}$*

(2.4) $\qquad \lambda \in \sigma(A_{BC}) \quad \Longleftrightarrow \quad 1 \in \sigma(\Delta_W(\lambda)),$

(2.5) $\qquad \lambda \in \sigma_*(A_{BC}) \quad \Longleftrightarrow \quad 1 \in \sigma_*(\Delta_W(\lambda)).$

(e) *If $1 \in \rho(\Delta_W(\lambda))$, then $\lambda \in \rho(A_{BC})$ and the resolvent of $A_{BC}$ is given by*

(2.6) $\qquad R(\lambda, A_{BC}) = R(\lambda, A) + R(\lambda, A^Z)B \cdot (Id_U - \Delta_U(\lambda))^{-1} \cdot CR(\lambda, A)$

(2.7) $\qquad = (Id_Z - \Delta_Z(\lambda))^{-1} \cdot R(\lambda, A).$

*Proof.* (a) follows from Lemma A.9 by choosing the spaces $F := Z_{-1}$, $E := X$ and operators $T := A_{BC}^Z$, $T_1 := T|_E = A_{BC}$. The assumption
$$E + \mathrm{rg}(T) = X + \mathrm{rg}(A_{BC}^Z) \quad \text{is dense in } F = Z_{-1},$$
which is needed for the last two inclusions follows by hypothesis. In fact, for $\nu \in \rho(A)$ the operator $(\nu - A^Z) : Z \to Z_{-1}$ is an isomorphism, hence the set
$$X + \mathrm{rg}(A_{BC}^Z) = X + (\nu - A_{BC}^Z)Z = (\nu - A^Z)\bigl(D(A) + (Id_Z - R(\nu, A^Z)BC)Z\bigr)$$
is dense in $Z_{-1}$ iff $D(A) + (Id_Z - R(\nu, A^Z)BC)Z$ is dense in $Z$.

We proceed by verifying (b)–(d) for $W = U$ and then return to the case $W = Z$ at the end.

To prove (b) we define for $\lambda \in \rho(A)$ the operator matrix
$$(2.8) \qquad \mathfrak{T} := \begin{pmatrix} \lambda - A^Z & B \\ C & Id_U \end{pmatrix} \in \mathcal{L}(Z \times U, Z_{-1} \times U).$$
Then the Schur complements of $\mathfrak{T}$ from Appendix A.4 are given by
$$\Delta_1 = \lambda - A_{BC}^Z \in \mathcal{L}(Z, Z_{-1}), \qquad \Delta_2 = Id_U - \Delta_U(\lambda) \in \mathcal{L}(U).$$
Hence, from Lemma A.10.(iv)–(vi) it follows that $\lambda - A_{BC}^Z$ is injective/has closed range/has dense range/has finite dimensional kernel/has range with finite co-dimension/is invertible iff $Id_U - \Delta_U(\lambda)$ has the same property, respectively. Since these properties characterize the various parts of the spectrum, this implies the first two equivalences in (b). The remaining ones then follow immediately from part (a).

For (c) assume that $\Delta_U(\lambda)$ is compact. Then using (a) and (b) we conclude
$$\lambda \in \sigma(A_{BC}) \implies \lambda \in \sigma(A_{BC}^Z) \iff 1 \in \sigma\bigl(\Delta_U(\lambda)\bigr) \iff 1 \in \sigma_p\bigl(\Delta_U(\lambda)\bigr)$$
$$\iff \lambda \in \sigma_p(A_{BC}^Z) \iff \lambda \in \sigma_p(A_{BC}) \implies \lambda \in \sigma(A_{BC}).$$
Therefore, all conditions are equivalent which proves the first chain of equivalences. If $U$ is finite dimensional, then $\Delta_U(\lambda)$ is compact for all $\lambda \in \rho(A)$ implying the second chain in (c).

For (d) assume that there exists $\nu \in \rho(A)$ such that $1 \in \rho(\Delta_U(\nu))$ which, by the first equivalence in (b), is equivalent to the existence of some $\nu \in \rho(A) \cap \rho(A_{BC}^Z)$. Then $A_{BC}^Z$ is closed and by Lemma 2.1 we conclude
$$(Z_{-1})_1^{A_{BC}^Z} \simeq Z \hookrightarrow X \hookrightarrow Z_{-1}.$$
Hence, for $F$, $E$, $T$ and $T_1$ as in the proof of (a) we have $(Z_{-1})_1^{A_{BC}^Z} = F_1^T \hookrightarrow E = X$. By Corollary A.9.(vii) this implies all equivalences for the various parts of the spectra.

To prove (e) assume that $1 \in \rho(\Delta_U(\lambda))$, i.e., $\lambda \in \rho(A_{BC}^Z)$. Then Lemma A.7.(vi) yields $R(\lambda, A_{BC}) = \Delta_1^{-1}|_X$. The formula for $\Delta_1^{-1}$ in Lemma A.10.(vi) then gives (2.6).

Finally, all assertions concerning $W = Z$ follow from Corollary A.12 applied to $E = U$, $F = Z$ and the operators $R := R(\lambda, A^Z)B \in \mathcal{L}(U, Z)$ and $Q := C \in \mathcal{L}(Z, U)$. □

*Remark* 2.4. (i) As main outcome, the previous result establishes that we always have
$$\lambda \in \sigma_p(A_{BC}) \iff 1 \in \sigma_p\bigl(\Delta_W(\lambda)\bigr)$$
while for the whole spectrum and its other parts only one implication holds in general. To obtain equivalence as for the point spectrum an additional assumption is necessary, e.g., that $Z_{-1} = X$, $\Delta_W(\lambda) \in \mathcal{L}(W)$ is compact or that (2.3) is satisfied. As Examples 2.5 and 3.11 show, such an extra assumption cannot be omitted in general.

(ii) If $U$ is finite dimensional, then (2.2) determines the spectral values of $A_{BC}$ in $\rho(A)$ as zeros of the *characteristic equation*
$$F(\lambda) := \det(Id_U - \Delta_U(\lambda)) = 0.$$
Since $F(\bullet)$ is holomorphic on $\rho(A)$, this implies that a connected component $\rho$ of $\rho(A)$ is either contained in $\sigma(A_{BC})$ (i.e. $F \equiv 0$ on $\rho$) or that $\sigma(A_{BC})$ has no accumulation in $\rho$ (i.e., $F \not\equiv 0$ on $\rho$). We note that even if $U$ is not finite dimensional but only a product space with a finite dimensional factor, using Schur complements might yield a characteristic equation as well. For concrete examples see Sections 4.3 and 4.4.

(iii) Note that by the previous result $A_{BC}$ is closed if $1 \in \rho(\Delta_W(\nu))$ for some $\nu \in \rho(A)$ (or, more generally, $1 \in \rho\big(\Delta_W(\nu_0,\nu)\big)$ for some $(\nu_0,\nu) \in \rho(A) \times \mathbb{C}$, cf. Corollary 2.7 below). This condition is in particular satisfied if $P = BC$ is a Weiss–Staffans perturbation of $A$, cf. [1, Def. 9], or if $\|CR(\nu, A^Z)B\| < 1$ for some $\nu \in \rho(A)$.

(iv) We mention that by definition of $A_{BC}$, Lemma A.7, Lemma A.10.(i)–(ii) and Corollary A.12 for $\lambda \in \rho(A)$ and $W \in \{U, Z\}$ we also have

$$\lambda - A_{BC}^Z \text{ surjective} \iff Id_W - \Delta_W(\lambda) \text{ surjective}$$
$$\implies \lambda - A_{BC} \text{ surjective.}$$

(v) We note that the operator $\mathcal{T}$ in (2.8) shows similarities to the system operator $S_\Sigma(\lambda)$ studied in some detail in [17].

The equivalences in (2.4) and (2.5) only hold with some extra assumption like (2.3) or compactness of $\Delta_W(\lambda)$. In fact, there are operators $A^Z$, $B$ and $C$ such that $A_{BC}^Z$ is not closed, hence $\sigma(A_{BC}^Z) = \mathbb{C}$, whereas $\sigma(A_{BC})$ might be rather small.

*Example* 2.5. For an unbounded, densely defined operator $(A, D(A))$ on a Banach space $X$ define on $\mathfrak{X} := X \times X$ the operator

$$\mathcal{G} := \begin{pmatrix} A & 0 \\ 0 & A \end{pmatrix} \cdot \begin{pmatrix} Id_X & Id_X \\ Id_X & Id_X \end{pmatrix}, \quad D(\mathcal{G}) := \left\{ \begin{pmatrix} x \\ y \end{pmatrix} \in X \times X : x + y \in D(A) \right\}.$$

Since the matrices

$$\begin{pmatrix} 1 & 1 \\ 1 & 1 \end{pmatrix} \quad \text{and} \quad \begin{pmatrix} 2 & 0 \\ 0 & 0 \end{pmatrix}$$

are similar, a simple computation shows that $\mathcal{G}$ is similar to

$$\mathcal{D} := \begin{pmatrix} 2A & 0 \\ 0 & 0 \end{pmatrix}, \quad D(\mathcal{D}) := D(A) \times X.$$

In particular, this implies that $\mathcal{G}$ is closed and $\sigma(\mathcal{G}) = \sigma(\mathcal{D}) = \sigma(2A) \cup \{0\}$.

A natural attempt to represent $\mathcal{G}$ as $\mathcal{G} = (\mathcal{A}^\mathcal{Z} + \mathcal{BC})|_\mathfrak{X}$ is to choose the spaces $\mathcal{Z} := \mathcal{U} := \mathfrak{X} = X \times X$, $\mathcal{Z}_{-1} := X_{-1} \times X_{-1}$ and the operators

$$\mathcal{A}^\mathcal{Z} := \begin{pmatrix} A_{-1} & 0 \\ 0 & A_{-1} \end{pmatrix} : \mathcal{Z} \subset \mathcal{Z}_{-1} \to \mathcal{Z}_{-1},$$

$$\mathcal{B} := \begin{pmatrix} 0 & A_{-1} \\ A_{-1} & 0 \end{pmatrix} \in \mathcal{L}(\mathcal{U}, \mathcal{Z}_{-1}) \quad \text{and} \quad \mathcal{C} := Id_\mathfrak{X} \in \mathcal{L}(\mathcal{Z}, \mathcal{U}).,$$

where $A_{-1} : X \subseteq X_{-1} \to X_{-1}$ denotes the extrapolated operator from Section A.2. Then a simple computation shows that $\rho(\mathcal{A}^\mathcal{Z}) = \rho(A)$ and $\mathcal{G} = (\mathcal{A}^\mathcal{Z} + \mathcal{BC})|_\mathfrak{X}$. However, neither there exists $\nu \in \rho(\mathcal{A}_{\mathcal{BC}}^\mathcal{Z})$, nor $1 \in \rho(\mathcal{C}R(\nu, \mathcal{A}^\mathcal{Z})\mathcal{B})$ for some $\nu \in \rho(\mathcal{A})$. In fact,

$$\mathcal{A}_{\mathcal{BC}}^\mathcal{Z} = \begin{pmatrix} A_{-1} & A_{-1} \\ A_{-1} & A_{-1} \end{pmatrix} : \mathcal{Z} \subset \mathcal{Z}_{-1} \to \mathcal{Z}_{-1}$$

is not closed, hence $\sigma(\mathcal{A}_{\mathcal{BC}}^\mathcal{Z}) = \mathbb{C}$. On the other hand, for $\lambda \in \rho(A)$

$$\mathcal{C}R(\lambda, \mathcal{A}^\mathcal{Z})\mathcal{B} = \begin{pmatrix} 0 & Id_X - \lambda R(\lambda, A) \\ Id_X - \lambda R(\lambda, A) & 0 \end{pmatrix} \in \mathcal{L}(\mathcal{U})$$

is never compact on $\mathfrak{X}$. Moreover, $Id_\mathcal{U} - \mathcal{C}R(\lambda, \mathcal{A}^\mathcal{Z})\mathcal{B}$ is invertible if and only if

$$1 \in \rho\big((\lambda R(\lambda, A) - Id_X)^2\big) \iff 0 \in \rho\big(\lambda^2 R(\lambda, A)^2 - 2\lambda R(\lambda, A)\big).$$

However, $\lambda^2 R(\lambda, A)^2 - 2\lambda R(\lambda, A)$ is never surjective, thus $1 \in \sigma(\mathcal{C}R(\lambda, \mathcal{A}^\mathcal{Z})\mathcal{B})$ as claimed.

Nevertheless, the matrix $\mathcal{G}$ can be treated also within our framework. To this end choose the spaces $\mathcal{Z} := \mathcal{U} := \mathfrak{X} = X \times X$ and

$$\mathcal{Z}_{-1} = \left\{ \begin{pmatrix} x \\ y \end{pmatrix} \in X_{-1} \times X_{-1} : x - y \in X \right\}$$

equipped with the norm induced by $X_{-1} \times X_{-1}$. Then we consider the operators

$$\mathcal{A}^{\mathcal{Z}} := \begin{pmatrix} 2A_{-1} & 0 \\ 2A_{-1} & 0 \end{pmatrix} : \mathcal{Z} \subset \mathcal{Z}_{-1} \to \mathcal{Z}_{-1},$$

$$\mathcal{B} := \begin{pmatrix} A_{-1} & 0 \\ A_{-1} & 0 \end{pmatrix} \in \mathcal{L}(\mathcal{U}, \mathcal{Z}_{-1}) \quad \text{and} \quad \mathcal{C} := \begin{pmatrix} -Id_X & Id_X \\ 0 & 0 \end{pmatrix} \in \mathcal{L}(\mathcal{Z}, \mathcal{U}).$$

For this choice we obtain $\rho(\mathcal{A}^{\mathcal{Z}}) = \rho(2A) \setminus \{0\}$ with resolvent

$$R(\lambda, \mathcal{A}^{\mathcal{Z}}) = \begin{pmatrix} R(\lambda, 2A) & 0 \\ R(\lambda, 2A) - \frac{1}{\lambda} & \frac{1}{\lambda} \end{pmatrix}, \quad \lambda \in \rho(\mathcal{A}^{\mathcal{Z}}).$$

Next one easily verifies that $\mathcal{G} = (\mathcal{A}^{\mathcal{Z}} + \mathcal{BC})|_{\mathcal{X}}$. Moreover, $Id_{\mathcal{U}} - \mathcal{C}R(\lambda, \mathcal{A}^{\mathcal{Z}})\mathcal{B} = Id_{\mathcal{U}}$ which is invertible for all $\lambda \in \rho(\mathcal{A}) = \rho(2A) \setminus \{0\}$. Hence, we can apply Theorem 2.3.(d) and conclude that $\sigma(\mathcal{G}) \subseteq \sigma(2A) \cup \{0\}$. We note that using Corollary 2.7 one could also show that $\sigma(\mathcal{G}) = \sigma(2A) \cup \{0\}$.

By applying Lemma A.10.(iii) to $\mathcal{T}$ in (2.8) we obtain using Corollary A.12 the following result generalizing [19, Thm. 1.1] and parts of [6, Lem. 6.4].

**Corollary 2.6.** *For $\lambda \in \rho(A)$ we have*

$$\ker(\lambda - A_{BC}) = R(\lambda, A^Z)B \ker(Id_U - \Delta_U(\lambda)) = \ker(Id_Z - \Delta_Z(\lambda)),$$
$$\ker(Id_U - \Delta_U(\lambda)) = C \ker(\lambda - A_{BC}).$$

One drawback of Theorem 2.3 is that it can be applied only to points $\lambda \in \rho(A)$. If one wants to determine the spectrum of a given operator $G$ it is therefore important to represent it as $G = A_{BC}$ for an operator $A$ having "small" spectrum. In many cases this is possible due to the great freedom in the choices of $A^Z$, $B$ and $C$ which only have to verify the Assumptions 1.1.

Nevertheless, we now present several approaches which allow to deal with points in $\lambda \in \sigma(A)$, too. The first one is based on the decomposition

$$(2.9) \qquad \lambda - A_{BC} = \lambda_0 - \left(A^Z + (BC + \lambda_0 - \lambda)\right)|_X$$

for some fixed $\lambda_0 \in \rho(A)$ and arbitrary $\lambda \in \mathbb{C}$. Define the "extended" space $\mathcal{U} := X \times U$ and for $\mu \in \mathbb{C}$ the "extended" operators

$$(2.10) \qquad \mathcal{B} := (Id_X, B) \in \mathcal{L}(\mathcal{U}, Z_{-1}), \qquad \mathcal{C}_\mu := \begin{pmatrix} \mu \cdot Id_Z \\ C \end{pmatrix} \in \mathcal{L}(Z, \mathcal{U}).$$

Then for $\mu := \lambda_0 - \lambda$ we obtain $\mathcal{BC}_\mu = BC + \lambda_0 - \lambda \in \mathcal{L}(Z, Z_{-1})$ which implies $A_{\mathcal{BC}_\mu} = (A^Z + (BC + \lambda_0 - \lambda))|_X$. Hence, by (2.9) we have $\lambda - A_{BC} = \lambda_0 - A_{\mathcal{BC}_\mu}$ and therefore

$$\lambda \in \sigma(A_{BC}) \iff \lambda_0 \in \sigma(A_{\mathcal{BC}_\mu}),$$
$$\lambda \in \sigma_*(A_{BC}) \iff \lambda_0 \in \sigma_*(A_{\mathcal{BC}_\mu})$$

for $* \in \{p, a, r, c, \text{ess}\}$, where now $\lambda_0 \in \rho(A)$. In order to apply our previous results to this situation we define the operators

$$\Delta_\mathcal{U}(\lambda_0, \lambda) := \mathcal{C}_\mu R(\lambda_0, A^Z)\mathcal{B} = \begin{pmatrix} (\lambda_0 - \lambda) \cdot R(\lambda_0, A) & (\lambda_0 - \lambda) \cdot R(\lambda_0, A^Z)B \\ CR(\lambda_0, A) & \Delta_U(\lambda_0) \end{pmatrix} \in \mathcal{L}(\mathcal{U}),$$

$$\Delta_Z(\lambda_0, \lambda) := R(\lambda_0, A^Z)\mathcal{BC}_\mu = \Delta_Z(\lambda_0) + (\lambda_0 - \lambda)R(\lambda_0, A) \in \mathcal{L}(Z)$$

for $\mu = \lambda_0 - \lambda$. Then by Theorem 2.3 the following holds.

**Corollary 2.7.** *Let $\lambda_0 \in \rho(A)$, $\lambda \in \mathbb{C}$ and $W \in \{\mathcal{U}, Z\}$.*

*(a) The following spectral implications always hold.*

$$\lambda \in \sigma(A_{BC}) \implies 1 \in \sigma(\Delta_W(\lambda_0, \lambda)),$$
$$\lambda \in \sigma_p(A_{BC}) \iff 1 \in \sigma_p(\Delta_W(\lambda_0, \lambda)),$$
$$\lambda \in \sigma_a(A_{BC}) \implies 1 \in \sigma_a(\Delta_W(\lambda_0, \lambda)),$$
$$\lambda \in \sigma_{ess}(A_{BC}) \implies 1 \in \sigma_{ess}(\Delta_W(\lambda_0, \lambda)).$$

If $D(A) + \text{rg}(Id_Z - \Delta_Z(\nu))$ is dense in $Z$ for some (hence all) $\nu \in \rho(A)$ then also

$$\lambda \in \sigma_c(A_{BC}) \quad \Longrightarrow \quad 1 \in \sigma_c(\Delta_W(\lambda_0, \lambda)),$$
$$\lambda \in \sigma_r(A_{BC}) \quad \Longleftarrow \quad 1 \in \sigma_r(\Delta_W(\lambda_0, \lambda)).$$

(b) If $\Delta_W(\lambda_0, \lambda) \in \mathcal{L}(W)$ is compact, then

$$\lambda \in \sigma(A_{BC}) \quad \Longleftrightarrow \quad \lambda \in \sigma_p(A_{BC}) \quad \Longleftrightarrow \quad 1 \in \sigma_p(\Delta_W(\lambda_0, \lambda)).$$

In particular, if $\dim(U) < \infty$, then

$$\lambda \in \sigma(A_{BC}) \quad \Longleftrightarrow \quad \lambda \in \sigma_p(A_{BC}) \quad \Longleftrightarrow \quad \det(Id_U - \Delta_U(\lambda_0, \lambda)) = 0.$$

(c) If there exist $\nu_0 \in \rho(A)$ and $\nu \in \mathbb{C}$ such that $1 \in \rho(\Delta_W(\nu_0, \nu))$, i.e., if $\rho(A_{BC}^Z) \neq \emptyset$, then

$$\lambda \in \sigma(A_{BC}) \quad \Longleftrightarrow \quad 1 \in \sigma(\Delta_W(\lambda_0, \lambda)),$$
$$\lambda \in \sigma_*(A_{BC}) \quad \Longleftrightarrow \quad 1 \in \sigma_*(\Delta_W(\lambda_0, \lambda))$$

for all $* \in \{a, r, c, ess\}$.

(d) If $1 \in \rho(\Delta_W(\lambda_0, \lambda))$, then $\lambda \in \rho(A_{BC})$ and the resolvent of $A_{BC}$ is given by

$$R(\lambda, A_{BC}) = R(\lambda_0, A) + R(\lambda_0, A^Z)\mathcal{B} \cdot (Id_\mathcal{U} - \Delta_\mathcal{U}(\lambda_0, \lambda))^{-1} \cdot \mathcal{C}_\mu R(\lambda_0, A)$$
$$= (Id_Z - \Delta_Z(\lambda_0, \lambda))^{-1} \cdot R(\lambda_0, A),$$

where $\mu := \lambda_0 - \lambda$.

In order to check the condition $1 \in \sigma(\Delta_\mathcal{U}(\lambda_0, \lambda))$ appearing above one might be tempted to use Schur complements, cf. Corollary A.12. To do so 1 has to be an element of the resolvent set of one of the diagonal entries of $\Delta_\mathcal{U}(\lambda_0, \lambda)$. However, the condition $1 \in \rho((\lambda_0 - \lambda) \cdot R(\lambda_0, A))$ is, by the spectral mapping theorem for the resolvent, see [7, IV.1.13], equivalent to $\lambda \in \rho(A)$. This case is already covered by Theorem 2.3. On the other hand, $1 \in \rho(\Delta_U(\lambda_0))$ is equivalent to $\lambda_0 \in \rho(A_{BC})$.

We now describe two other approaches which allow to deal with boundary points $\lambda_0 \in \partial\sigma(A)$. This is in particular useful if $\sigma(A) = \partial\sigma(A)$, e.g., if $\sigma(A)$ is discrete in $\mathbb{C}$, e.g., if $A$ has compact resolvent. But also in other cases this might be very helpful since $\partial\sigma(A)$ already determines the *spectral bound*

$$s(A) := \sup\{\operatorname{Re}\lambda : \lambda \in \sigma(A)\}$$

of an operator $A$ which is closely related to the asymptotic behavior of the solutions of the abstract Cauchy problem (ACP), cf. [7, Cor. IV.3.12].

**Corollary 2.8.** *Let $\lambda_0 \in \partial\sigma(A)$. If $\Delta_Z(\bullet) = R(\bullet, A^Z)BC : \rho(A) \to \mathcal{L}(Z)$ has a continuous extension $\bar\Delta_Z(\bullet)$ in $\lambda_0$, then $\lambda_0 \in \sigma(A_{BC})$.*

*Proof.* By contradiction assume that $\lambda_0 \in \rho(A_{BC})$. Define

$$R(\bullet) : (\rho(A) \cap \rho(A_{BC})) \cup \{\lambda_0\} \to \mathcal{L}(X), \quad R(\lambda) := (Id_Z - \bar\Delta_Z(\lambda)) \cdot R(\lambda, A_{BC}).$$

Then by (2.7) we conclude $R(\lambda) = R(\lambda, A)$ for all $\lambda \in \rho := \rho(A) \cap \rho(A_{BC})$, i.e., for all $\lambda \in \rho$ we have

$$R(\lambda)(\lambda - A)x = x \quad \text{for all } x \in D(A),$$
$$(\lambda - A^Z)R(\lambda)x = x \quad \text{for all } x \in X.$$

Since $\lambda_0 \in \bar\rho$, by continuity these relations remain valid for $\lambda = \lambda_0$. This implies $R(\lambda_0)x \in D(A)$ and therefore $R(\lambda_0) = R(\lambda_0, A) \in \mathcal{L}(X)$. This proves $\lambda_0 \in \rho(A)$ which contradicts the assumption that $\lambda \in \partial\sigma(A) \subset \sigma(A)$. Hence, $\lambda_0 \in \sigma(A_{BC})$ as claimed. $\square$

The following result is a "local" version of the previous one. In fact, we do not suppose that $\Delta_Z(\bullet)$ has a continuous extension to some spectral value on the boundary of $\sigma(A)$ but only $\Delta_Z(\bullet)f$ for an associated eigenvector $f$ of $A$.

**Corollary 2.9.** *Let $\lambda_0 \in \partial\sigma(A) \cap \sigma_p(A)$. If for some $0 \neq f \in \ker(\lambda_0 - A)$ the limit*

$$\lim_{\rho(A) \ni \lambda \to \lambda_0} \Delta_Z(\lambda) f =: g$$

*converges in $X$ such that $g \in Z$ and $Cg = 0$, then $\lambda_0 \in \sigma_p(A_{BC})$.*

*Proof.* We first show that $f + g \in \ker(\lambda_0 - A_{BC})$. Indeed, $f + g \in Z = D(A_{BC}^Z)$ and

$$\begin{aligned}
(\lambda_0 - A_{BC}^Z)(f+g) &= (\lambda_0 - A^Z)g - BCf \\
&= (\lambda_0 - A^Z)(g - \Delta_Z(\lambda)f) + (\lambda_0 - \lambda) \cdot \Delta_Z(\lambda)f \\
&\to 0 \quad \text{in } Z_{-1} \text{ as } \rho(A) \ni \lambda \to \lambda_0,
\end{aligned}$$
(2.11)

where we used that $A^Z \in \mathcal{L}(X, Z_{-1})$. This implies $f + g \in \ker(\lambda_0 - A_{BC}^Z) = \ker(\lambda_0 - A_{BC})$ as claimed. Now, if $f + g \neq 0$ it follows immediately that $\lambda_0 \in \sigma_p(A_{BC})$. On the other hand, if $f + g = 0$, i.e., $f = -g$, then by the assumptions $Cg = 0$ and $0 \neq f \in \ker(\lambda_0 - A)$ we obtain

$$(\lambda_0 - A_{BC})f = (\lambda_0 - A^Z - BC)f = 0.$$

This implies again $\lambda_0 \in \sigma_p(A_{BC})$, hence the proof is complete. □

## 3. The Generic Example

In this section we introduce a general setting which generalizes "boundary" perturbations of operators in the sense of Greiner, cf. [8], and then apply to it the theory developed in Section 2. Concrete application fitting into this framework can be found in Section 4.

**3.1. The Operator $A_P^\Phi$.** We start with a Banach space $X$ and a linear "*maximal operator*" (maximal in the sense of a "big" domain, e.g., a differential operator without boundary conditions) $A_m : D(A_m) \subseteq X \to X$. In order to single out a restriction $A$ of $A_m$ we take a Banach spaces $\partial X$, called "*space of boundary conditions*", and a linear "*boundary operator*", e.g. a "trace" operator, $L : D(A_m) \to \partial X$ and define on $X$

(3.1) $\qquad A \subseteq A_m, \quad D(A) = \{f \in D(A_m) : Lf = 0\} = \ker(L).$

Next we perturb $A$ in the following way. For operators $P : D(A_m) \to X$ and $\Phi : D(A_m) \to \partial X$ we introduce $A_P^\Phi : D(A_P^\Phi) \subseteq X \to X$ given by

(3.2) $\qquad A_P^\Phi \subseteq A_m + P, \quad D(A_P^\Phi) := \{f \in D(A_m) : Lf = \Phi f\} = \ker(L - \Phi).$

Hence, $A_P^\Phi$ can be considered as a twofold perturbation of $A$,

- by the operator $P$ to change its action, and
- by the operator $\Phi$ to change its domain,

cf. Diagram 2. We note that in [8] the operator $\Phi : X \to \partial X$ has to be bounded and $P = 0$.

$$\begin{CD}
X \supseteq D(A_m) @>{A_m,\, P}>> X \\
@. @VV{L,\, \Phi}V \\
@. \partial X
\end{CD}$$

Diagram 2: The operators defining $A_P^\Phi$ in (3.2).

Under the Assumptions 3.5 below, which cover unbounded $\Phi$ and $P$, the spectral properties of $A_P^\Phi$ can be studied using our results from Section 2. As a first step towards this goal we introduce in the next subsection so-called "abstract Dirichlet operators" and then study their existence and basic properties.

3.2. **Abstract Dirichlet Operators.** If for some $\lambda \in \mathbb{C}$ the restriction
$$L|_{\ker(\lambda - A_m)} : \ker(\lambda - A_m) \to \partial X$$
is invertible with inverse
$$L_\lambda := \left(L|_{\ker(\lambda - A_m)}\right)^{-1} : \partial X \to X,$$
then we call $L_\lambda$ the abstract *Dirichlet operator* associated to $\lambda$, $A_m$ and $L$. This notion is motivated by the fact that for a given "boundary value" $x \in \partial X$ the "function" $f = L_\lambda x$ is the unique solution of the "abstract Dirichlet Problem"

(aDP) $$\begin{cases} (\lambda - A_m)f = 0, \\ Lf = x. \end{cases}$$

Our approach is mainly based on these Dirichlet operators $L_\lambda$ and in general we do not have to know the resolvent $R(\lambda, A)$ explicitly. Therefore, the following result characterizing resolvent points of $A$ in terms of the existence of $L_\lambda$ might be helpful in applications.

**Lemma 3.1.** *Let $A$ be given by (3.1) and assume that $L : D(A_m) \to \partial X$ is surjective. Then $\lambda \in \rho(A)$ if and only if*
 (i) *$A$ is closed,*
 (ii) *$\lambda - A_m : D(A_m) \to X$ is surjective,*
 (iii) *$L_\lambda : \partial X \to X$ exists (algebraically), i.e., for every $x \in \partial X$ the abstract Dirichlet Problem (aDP) has a unique solution $f = L_\lambda x \in D(A_m)$.*

*Proof.* If $\lambda \in \rho(A)$ then clearly $A$ is closed and $\lambda - A_m$ is surjective. To show (iii) fix some $x \in \partial X$. Then there exists $h \in D(A_m)$ such that $Lh = x$. Define
$$f := \bigl(Id_X - R(\lambda, A)(\lambda - A_m)\bigr)h \in D(A_m).$$
Since $D(A) = \ker(L)$, we have $Lf = x$. Moreover, from $(\lambda - A_m)R(\lambda, A) = Id_X$ we conclude $f \in \ker(\lambda - A_m)$, i.e., $f$ solves (aDP). Now assume that $f_1, f_2 \in D(A_m)$ are both solutions of (aDP). Then $(\lambda - A_m)(f_1 - f_2) = 0$ and $L(f_1 - f_0) = 0$, i.e., $f_1 - f_2 \in \ker(\lambda - A) = \{0\}$ which shows uniqueness.

Conversely, assume that (i)–(iii) hold. First we show that $\lambda - A$ is surjective. Let $g \in X$, then by (ii) there exist $h \in D(A_m)$ such that $(\lambda - A_m)h = g$. Let $f := (Id_X - L_\lambda L)h$. Then $Lf = 0$, i.e., $f \in D(A)$ and
$$(\lambda - A)f = (\lambda - A_m)(Id_X - L_\lambda L)h = g$$
showing surjectivity. To prove injectivity let $(\lambda - A)f = 0$ for some $f \in D(A) = \ker(L)$. Then $f \in \ker(\lambda - A_m)$ and $Lf = 0$. Since the same holds for $f = 0$, by the uniqueness assumption in (iii) we conclude $f = 0$. Summing up, this shows that $\lambda - A$ is bijective and since $A$ by (i) is closed, the closed graph theorem implies $\lambda \in \rho(A)$ as claimed. $\square$

Next we give a closedness condition ensuring existence and boundedness of the Dirichlet operators. More precisely, let

(3.3) $$\mathcal{A} := \binom{A_m}{L} : D(A_m) \subset X \to X \times \partial X$$

and define $Z := (D(A_m), \|\cdot\|_\mathcal{A}) = [D(\mathcal{A})]$. Then the following holds.

**Lemma 3.2.** *Assume that $L : D(A_m) \to \partial X$ is surjective and $\rho(A) \neq \emptyset$. Then the following conditions are equivalent.*
 (a) *The operator $\mathcal{A}$ in (3.3) is closed, i.e., $Z$ is a Banach space.*
 (b) *The Dirichlet operator $L_\lambda \in \mathcal{L}(\partial X, X)$ exists for all $\lambda \in \rho(A)$.*
 (c) *The Dirichlet operator $L_\lambda \in \mathcal{L}(\partial X, X)$ exists for some $\lambda \in \rho(A)$.*

*Proof.* (a)$\Rightarrow$(b) is shown in [5, Lem. 2.2] while (b)$\Rightarrow$(c) is trivially satisfied. To prove (c)$\Rightarrow$(a) take a sequence $(f_n)_{n \in \mathbb{N}} \subset D(A_m)$ such that $f_n \to f_0 \in X$ and
$$\binom{A_m}{L}f_n \to \binom{g_0}{x_0} \in X \times \partial X \qquad \text{as } n \to +\infty.$$

Since $L_\lambda$ is bounded, this implies $D(A) = \ker(L) \ni (Id - L_\lambda L)f_n \to f_0 - L_\lambda x_0$ and
$$(\lambda - A)(Id - L_\lambda L)f_n = (\lambda - A_m)f_n \to \lambda f_0 - g_0.$$
Hence, by the closedness of $\lambda - A$ we conclude $f_0 - L_\lambda x_0 \in D(A)$ and $(\lambda - A)(f_0 - L_\lambda x_0) = \lambda f_0 - g_0$. Using this we obtain $f_0 \in \ker(\lambda - A_m) + D(A) \subseteq D(A_m)$ and $(\lambda - A_m)f_0 = \lambda f_0 - g_0$, i.e., $A_m f_0 = g_0$. Finally, from $f_0 - L_\lambda x_0 \in D(A) = \ker(L)$ it follows that $Lf_0 = x_0$. Summing up this proves (a). □

*Remark* 3.3. Note that if $A_m : D(A_m) \subseteq X \to X$ is closed and $L \in \mathcal{L}([D(A_m)], \partial X)$ is bounded then $\binom{A_m}{L} : D(A_m) \to X \times \partial X$ is closed. Hence, by the previous result, closedness of $A_m$ and relative boundedness of $L$ imply the existence of the Dirichlet operators $L_\lambda$ for all $\lambda \in \rho(A)$.

We need two more facts concerning the Dirichlet operators which follow as in [2, Prop. 3.2] and [8, Lem. 1.2], respectively.

**Lemma 3.4.** *Assume that for some $\lambda \in \mathbb{C}$ the Dirichlet operator $L_\lambda : D(A_m) \to \partial X$ exists. Then*
 (i) *for every $\mu \in \rho(A)$ also $L_\mu : D(A_m) \to \partial X$ exists and is given by*
(3.4) $$L_\mu = (\lambda - A)R(\mu, A)L_\lambda;$$
 (ii) *the domain of $A_m$ can be decomposed as algebraically direct sum*
$$D(A_m) = D(A) \oplus \ker(\lambda - A)$$
 *with projection $P_\lambda := L_\lambda L$ satisfying $\operatorname{rg}(P_\lambda) = \ker(\lambda - A)$ and $\ker(P_\lambda) = D(A)$.*

3.3. **How to represent $A_P^\Phi$ as $A_{BC}$?** In order to apply the results of Section 2 to the generic example we have to represent $A_P^\Phi$ from (3.2) as $(A^Z + BC)|_X$ like in (1.1) for suitable operators $A^Z$, $B$ and $C$ satisfying Assumptions 1.1. In order to do so we make in the remaining part of this section the following standing

*Assumptions* 3.5. (i) $A : D(A) \subseteq X \to X$ satisfies $\rho(A) \neq \emptyset$.
 (ii) $L : D(A_m) \to \partial X$ is surjective.
 (iii) One of the (equivalent) conditions (a)–(c) in Lemma 3.2 is satisfied.
 (iv) $\Phi \in \mathcal{L}(Z, \partial X)$ and $P \in \mathcal{L}(Z, X)$ for the Banach space $Z := [D(\mathcal{A})]$, cf. (3.3).

Note that assumption (iv) is equivalent to the existence of $M \geq 0$ such that
$$\|\Phi f\|_{\partial X} \leq M \cdot (\|A_m f\|_X + \|Lf\|_{\partial X} + \|f\|_X),$$
$$\|Pf\|_X \leq M \cdot (\|A_m f\|_X + \|Lf\|_{\partial X} + \|f\|_X)$$
for all $f \in D(A_m)$.

Our first aim is now to extend the operator $A : D(A) \subseteq X \to X$ to $A^Z : Z \subseteq Z_{-1} \to Z_{-1}$ for $Z = D(A_m)$ and a space $Z_{-1}$ containing $X$.

Consider on $\ker(\lambda - A_m)$ the norm
$$\|f\|_L := \|f\|_X + \|Lf\|_{\partial X}.$$
Then the following holds.

**Lemma 3.6.** *For every $\lambda \in \mathbb{C}$ the space $(\ker(\lambda - A_m), \|\bullet\|_L)$ is complete.*

*Proof.* By Assumption 3.5.(iii) the operator $\mathcal{A} := \binom{A_m}{L}$ is closed, hence $Z := (D(A_m), \|\bullet\|_\mathcal{A})$ is complete. Moreover, $\lambda - A_m : Z \to X$ is bounded, hence $\ker(\lambda - A_m)$ is closed in $Z$. However, on $\ker(\lambda - A_m)$ the norms $\|\bullet\|_\mathcal{A}$ and $\|\bullet\|_L$ are equivalent which implies the claim. □

We proceed and define for some fixed $\eta_0 \in \rho(A)$ the Banach space
$$Z_{-1} := X \times \ker(\eta_0 - A_m) \qquad \text{equipped with the norm} \qquad \|\binom{f}{k}\|_{Z_{-1}} := \|f\|_X + \|k\|_L.$$
Identifying $X$ by $X \times \{0\}$, i.e., $f \in X$ by $(f, 0)^\top \in Z_{-1}$, we obtain $X \hookrightarrow Z_{-1}$. Recall that by Lemma 3.4.(ii) we have $Z = D(A_m) = D(A) \oplus \ker(\eta_0 - A_m)$ with corresponding projection $P_{\eta_0} = L_{\eta_0} L \in \mathcal{L}(Z)$. Using this we define the extension $A^Z : Z \subseteq Z_{-1} \to Z_{-1}$ by
$$A^Z f := \binom{A(Id - P_{\eta_0})f}{P_{\eta_0} f}$$
cf. Diagram 3.

$$\begin{array}{ccc} Z_{-1} = & X & \times \ker(\eta_0 - A_m) \\ A^Z \uparrow & A \uparrow & Id \uparrow \\ Z = & D(A) & \oplus \ker(\eta_0 - A_m) \end{array}$$

Diagram 3: The extension $A^Z$ of $A$ to $Z$.

**Lemma 3.7.** *We have $A = A^Z|_X$ and $R(\lambda, A) = R(\lambda, A^Z)|_X$ for all $\lambda \in \rho(A) = \rho(A^Z)$.*

*Proof.* The first assertion follows from the fact that $P_{\eta_0} f = 0$ if and only if $f \in D(A)$. To show the remaining assertions we first prove that $\eta_0 \in \rho(A^Z)$. To this end define $R : Z_{-1} \to Z_{-1}$ by

$$R\binom{f}{k} := R(\eta_0, A)(f + \eta_0 k) - k \in Z = D(A^Z).$$

Using that $P_{\eta_0}|_{\ker(\eta_0 - A_m)} = Id$ and $P_{\eta_0} R(\eta_0, A) = 0$ we obtain

$$(\eta_0 - A^Z) \cdot R\binom{f}{k} = \binom{(\eta_0 - A(Id - P_{\eta_0}))\big(R(\eta_0, A)(f + \eta_0 k) - k\big)}{-P_{\eta_0}\big(R(\eta_0, A)(f + \eta_0 k) - k\big)} = \binom{f}{k} \quad \text{for all } \binom{f}{k} \in Z_{-1},$$

$$R \cdot (\eta_0 - A^Z) f = R(\eta_0, A) \cdot \big(\eta_0 f - A(Id - P_{\eta_0}) f - \eta_0 P_{\eta_0} f\big) + P_{\eta_0} f = f \quad \text{for all } f \in Z,$$

i.e., $R = (\eta_0 - A^Z)^{-1}$ is the algebraic inverse. Moreover, $R \in \mathcal{L}(Z_{-1})$ and hence $\eta_0 \in \rho(A^Z) \neq \emptyset$. Since $Z = D(A^Z) \hookrightarrow X$, Lemma A.7.(vi)–(vii) applied to $E = X$, $F = Z_{-1}$ and $T = A^Z$ then implies $\rho(A^Z) = \rho((A^Z)_1) = \rho(A)$ and $R(\lambda, A) = R(\lambda, A^Z)|_X$ as claimed. □

To proceed we define the operator

(3.5) $$L_A := (\lambda - A^Z) L_\lambda \in \mathcal{L}(\partial X, Z_{-1})$$

which by (3.4) is independent of $\lambda \in \rho(A)$. Moreover, we introduce $U := X \times \partial X$ and

(3.6) $$B := (Id_X, L_A) \in \mathcal{L}(U, Z_{-1}) \quad \text{and} \quad C_\mu := \binom{P + \mu \cdot Id_Z}{\Phi} \in \mathcal{L}(Z, U)$$

for $\mu \in \mathbb{C}$. Then we obtain the following representation of $A_P^\Phi$, cf. Diagram 4.

**Lemma 3.8.** *Let $\lambda_0, \lambda \in \mathbb{C}$ and $\mu := \lambda_0 - \lambda$. Then for $A_P^\Phi$ given by (3.2) one has*

(3.7) $$\lambda - A_P^\Phi = \lambda_0 - \big(A^Z + P + (\lambda_0 - \lambda) + L_A \Phi\big)\big|_X = \lambda_0 - A_{BC_\mu}.$$

Diagram 4: The operators appearing in the representation $A_P^\Phi = (A^Z + P + L_A \Phi)|_X$.

*Proof.* To prove (3.7) it suffices to consider the case $\lambda = \lambda_0 = \mu = 0$, i.e., to verify that

(3.8) $$A_P^\Phi = \big(A^Z + P + L_A \Phi\big)\big|_X.$$

Denote by $G$ the operator defined by the right-hand-side of (3.8) and fix some $\lambda \in \rho(A)$. Then for $f \in D(A_m)$ we have

$$\begin{aligned} f \in D(G) &\iff (A^Z - \lambda)(Id - L_\lambda \Phi)f + (P + \lambda)f \in X \\ &\iff (Id - L_\lambda \Phi)f \in D(A) = \ker(L) \\ &\iff Lf = \Phi f \\ &\iff f \in D(A_P^\Phi). \end{aligned}$$

Moreover, for $f \in D(A_P^\Phi)$ we obtain

$$\begin{aligned} Gf &= (A - \lambda)(Id - L_\lambda \Phi)f + (P + \lambda)f \\ &= (A_m - \lambda)f + (P + \lambda)f \\ &= (A_m + P)f = A_P^\Phi f, \end{aligned}$$

hence $A_P^\Phi = G$ as claimed. $\square$

Under suitable assumptions the spectra of $A^Z + P + L_A \Phi : D(A_m) \subseteq Z_{-1} \to Z_{-1}$ and its restriction $A_P^\Phi$ to $X$ coincide. More precisely we have the following which generalizes Lemma 3.7.

**Lemma 3.9.** *Assume that $L - \Phi : D(A_m) \to \partial X$ is surjective and that*

$$\begin{pmatrix} A_m + P \\ L - \Phi \end{pmatrix} : D(A_m) \subset X \to X \times \partial X$$

*is closed. Then*

$$\rho(A_P^\Phi) = \rho(A^Z + P + L_A \Phi).$$

*Proof.* The inclusion "$\supseteq$" is clear by Lemma A.7.(vi). Hence we only have to verify that $\lambda \in \rho(A_P^\Phi)$ implies $\lambda \in \rho(A^Z + P + L_A \Phi)$. To this end note that the pair $A_m + P : D(A_m) \subseteq X \to X$, $L - \Phi : D(A_m) \to \partial X$ satisfies assumption (a) of Lemma 3.2. Hence, the corresponding Dirichlet operator

$$K_\lambda := \left((L - \Phi)|_{\ker(\lambda - A_m - P)}\right)^{-1} \in \mathcal{L}(\partial X, X)$$

exists. Using this we define $R_\lambda \in \mathcal{L}(Z_{-1})$ by

$$R_\lambda \begin{pmatrix} f \\ k \end{pmatrix} := R(\lambda, A_P^\Phi)(f + \eta_0 k) - K_\lambda Lk \in Z = D(A^Z + P + L_A \Phi).$$

Then a simple computation as in the proof of Lemma 3.7 shows that $R_\lambda = R(\lambda, A^Z + P + L_A \Phi)$, i.e., $\lambda \in \rho(A^Z + P + L_A \Phi)$ as claimed. $\square$

3.4. **Spectral Theory for $A_P^\Phi$.** By applying Theorem 2.3 to the representation (3.7) of $\lambda - A_P^\Phi$ we easily obtain the following result where part (c) follows by using also Lemma 3.9.

**Corollary 3.10.** *For $\lambda_0 \in \rho(A)$, $\lambda \in \mathbb{C}$ and $\mu := \lambda_0 - \lambda$ define on $\mathcal{U} := X \times \partial X$ and $Z$ the operators*

$$\Delta_\mathcal{U}(\lambda_0, \lambda) := \mathcal{C}_\mu R(\lambda_0, A^Z)\mathcal{B} = \begin{pmatrix} (P + \mu \cdot Id_Z)R(\lambda_0, A) & (P + \mu \cdot Id_Z)L_{\lambda_0} \\ \Phi R(\lambda_0, A) & \Phi L_{\lambda_0} \end{pmatrix} \in \mathcal{L}(\mathcal{U}),$$

$$\Delta_Z(\lambda_0, \lambda) := R(\lambda_0, A^Z)\mathcal{B}\mathcal{C}_\mu = R(\lambda_0, A)(P + \mu \cdot Id_Z) + L_{\lambda_0}\Phi \in \mathcal{L}(Z).$$

*Then for $W \in \{\mathcal{U}, Z\}$ the following holds.*

*(a) We always have*

$$\begin{aligned} \lambda \in \sigma(A_P^\Phi) &\implies 1 \in \sigma(\Delta_W(\lambda_0, \lambda)), \\ \lambda \in \sigma_p(A_P^\Phi) &\iff 1 \in \sigma_p(\Delta_W(\lambda_0, \lambda)), \\ \lambda \in \sigma_a(A_P^\Phi) &\implies 1 \in \sigma_a(\Delta_W(\lambda_0, \lambda)), \\ \lambda \in \sigma_{ess}(A_P^\Phi) &\implies 1 \in \sigma_{ess}(\Delta_W(\lambda_0, \lambda)). \end{aligned}$$

If $D(A) + \mathrm{rg}(Id_Z - L_\nu \Phi)$ is dense in $Z$ for some (hence all) $\nu \in \rho(A)$ then also
$$\lambda \in \sigma_c(A_P^\Phi) \quad \Longrightarrow \quad 1 \in \sigma_c\big(\Delta_W(\lambda_0, \lambda)\big),$$
$$\lambda \in \sigma_r(A_P^\Phi) \quad \Longleftarrow \quad 1 \in \sigma_r\big(\Delta_W(\lambda_0, \lambda)\big).$$

(b) If $\Delta_W(\lambda_0, \lambda) \in \mathcal{L}(W)$ is compact, then
$$\lambda \in \sigma(A_P^\Phi) \iff \lambda \in \sigma_p(A_P^\Phi) \iff 1 \in \sigma_p\big(\Delta_W(\lambda_0, \lambda)\big).$$

(c) If $L - \Phi : D(A_m) \to \partial X$ is surjective and $\binom{A_m + P}{L - \Phi} : D(A_m) \subset X \to X \times \partial X$ is closed, then

(3.9) $$\lambda \in \sigma(A_P^\Phi) \iff 1 \in \sigma\big(\Delta_W(\lambda_0, \lambda)\big).$$

(d) If there exist $\nu_0 \in \rho(A)$, $\nu \in \mathbb{C}$ such that $1 \in \rho\big(\Delta_W(\nu_0, \nu)\big)$, then for all $* \in \{a, r, c, ess\}$
$$\lambda \in \sigma(A_P^\Phi) \iff 1 \in \sigma\big(\Delta_W(\lambda_0, \lambda)\big),$$
$$\lambda \in \sigma_*(A_P^\Phi) \iff 1 \in \sigma_*\big(\Delta_W(\lambda_0, \lambda)\big).$$

(e) If $1 \in \rho\big(\Delta_W(\lambda_0, \lambda)\big)$, then $\lambda \in \rho(A_P^\Phi)$ and the resolvent of $A_P^\Phi$ is given by
$$R(\lambda, A_P^\Phi) = R(\lambda_0, A) + \big(R(\lambda_0, A), L_{\lambda_0}\big) \cdot \big(Id_\mathcal{U} - \Delta_\mathcal{U}(\lambda_0, \lambda)\big)^{-1} \cdot \binom{(P + \mu \cdot Id_Z)R(\lambda_0, A)}{\Phi R(\lambda_0, A)}$$
$$= \big(Id_Z - R(\lambda_0, A)(P + \mu \cdot Id_Z) - L_{\lambda_0}\Phi\big)^{-1} \cdot R(\lambda_0, A),$$
where $\mu = \lambda_0 - \lambda$.

The following simple example shows that the equivalence in part (c) does not hold without the surjectivity assumption on $L - \Phi$.

*Example* 3.11. On $X := \mathrm{L}^p(\mathbb{R}_+)$ consider $A_m := \frac{d}{ds}$ with domain $D(A_m) := \mathrm{W}^{1,p}(\mathbb{R}_+)$. Moreover, take $\partial X := \mathbb{C}$ and $L := \delta_0 : D(A_m) \to \partial X$. Then $A_m$ is closed, $L$ is surjective and for $A := A_m|_{\ker(L)}$ we have $\sigma(A) = \{\zeta \in \mathbb{C} : \mathrm{Re}\,\zeta \geq 0\}$. Next, choose $P = 0$ and $\Phi = \delta_0 = L$ which gives $A_0^\Phi =: A^\Phi = A_m|_{\ker(L-\Phi)} = A_m$ and $\sigma(A^\Phi) = \{\zeta \in \mathbb{C} : \mathrm{Re}\,\zeta \leq 0\}$.

In this case we obtain for arbitrary $\lambda_0 \in \rho(A)$ and $\lambda \in \mathbb{C}$
$$\Delta_\mathcal{U}(\lambda_0, \lambda) = \begin{pmatrix} (\lambda_0 - \lambda)R(\lambda_0, A) & (\lambda_0 - \lambda)L_{\lambda_0} \\ 0 & 1 \end{pmatrix}.$$

This implies that $1 \in \sigma(\Delta_\mathcal{U}(\lambda_0, \lambda))$ for all $\lambda_0 \in \rho(A)$, $\lambda \in \mathbb{C}$ while $\lambda \in \sigma(A^\Phi)$ iff $\mathrm{Re}\,\lambda \leq 0$. However, in this case $L - \Phi = 0$ is not surjective and hence these facts do not contradict part (c) of the previous result.

Next we give a simple condition ensuring that $L - \Phi$ is surjective. For an application in the context of delay equations see Corollary 4.11.

**Lemma 3.12.** *If there exists $\lambda_0 \in \rho(A)$ such that*

(3.10) $$\Phi\big(\ker(\lambda_0 - A_m)\big) \subseteq \Phi\big(D(A)\big),$$

*then $L - \Phi : D(A_m) \to \partial X$ is surjective.*

*Proof.* It suffices to show that for each $x \in \partial X$ there exists $f \in D(A_m)$ such that

(3.11) $$\begin{cases} Lf = x, \\ \Phi f = 0. \end{cases}$$

Let $x \in \partial X$. Then by (3.10) there exists $f_0 \in D(A) = \ker(L)$ such that $\Phi(L_{\lambda_0}x + f_0) = 0$. This implies that $f := L_{\lambda_0}x + f_0$ solves (3.11) and the proof is complete. $\square$

By choosing $\Phi = \alpha L$ for some $\alpha \neq 1$ it is clear that (3.11) is only sufficient but not necessary for $L - \Phi$ to be surjective.

If $P = 0$, we obtain the operator $A^\Phi := A_0^\Phi \subset A_m$ with domain $D(A^\Phi) = \ker(L - \Phi)$. If also $\lambda = \lambda_0 \in \rho(A)$ we can cancel out the unnecessary terms and consider $U = \partial X$, $B = L_A$ and $C = \Phi$. Then $A^\Phi = A_{BC}$ and the previous result simplifies as follows.

**Corollary 3.13.** *For $\lambda \in \rho(A)$ define the operators*
$$\Delta_{\partial X}(\lambda) = CR(\lambda, A^Z)B = \Phi L_\lambda \in \mathcal{L}(\partial X),$$
$$\Delta_Z(\lambda) = R(\lambda, A^Z)BC = L_\lambda \Phi \in \mathcal{L}(Z).$$
*Then for $W \in \{\partial X, Z\}$ the following holds.*

(a) *We always have*
$$\begin{aligned}
\lambda \in \sigma(A^\Phi) &\implies& 1 \in \sigma\big(\Delta_W(\lambda)\big), \\
\lambda \in \sigma_p(A^\Phi) &\iff& 1 \in \sigma_p\big(\Delta_W(\lambda)\big), \\
\lambda \in \sigma_a(A^\Phi) &\implies& 1 \in \sigma_a\big(\Delta_W(\lambda)\big), \\
\lambda \in \sigma_{ess}(A^\Phi) &\implies& 1 \in \sigma_{ess}\big(\Delta_W(\lambda)\big).
\end{aligned}$$
*If $D(A) + \mathrm{rg}(Id_Z - L_\nu \Phi)$ is dense in $Z$ for some (hence all) $\nu \in \rho(A)$ then also*
$$\begin{aligned}
\lambda \in \sigma_c(A^\Phi) &\implies& 1 \in \sigma_c\big(\Delta_W(\lambda)\big), \\
\lambda \in \sigma_r(A^\Phi) &\impliedby& 1 \in \sigma_r\big(\Delta_W(\lambda)\big).
\end{aligned}$$

(b) *If $\Delta_W(\lambda) \in \mathcal{L}(W)$ is compact, then*

(3.12) $\qquad \lambda \in \sigma(A^\Phi) \iff \lambda \in \sigma_p(A^\Phi) \iff 1 \in \sigma_p\big(\Delta_W(\lambda)\big).$

*In particular, if $\dim(\partial X) < \infty$, then*

(3.13) $\qquad \lambda \in \sigma(A^\Phi) \iff \lambda \in \sigma_p(A^\Phi) \iff \det\big(Id_{\partial X} - \Delta_{\partial X}(\lambda)\big) = 0.$

(c) *If $L - \Phi : D(A_m) \to \partial X$ is surjective and $\binom{A_m}{L-\Phi} : D(A_m) \subset X \to X \times \partial X$ is closed, then*
$$\lambda \in \sigma(A^\Phi) \iff 1 \in \sigma\big(\Delta_W(\lambda)\big).$$

(d) *If there exist $\nu_0 \in \rho(A)$ such that $1 \in \rho\big(\Delta_W(\nu_0)\big)$, then for all $* \in \{a, r, c, ess\}$*
$$\begin{aligned}
\lambda \in \sigma(A^\Phi) &\iff& 1 \in \sigma\big(\Delta_W(\lambda)\big), \\
\lambda \in \sigma_*(A^\Phi) &\iff& 1 \in \sigma_*\big(\Delta_W(\lambda)\big).
\end{aligned}$$

(e) *If $1 \in \rho\big(\Delta_W(\lambda)\big)$, then $\lambda \in \rho(A^\Phi)$ and the resolvent of $A^\Phi$ is given by*
$$\begin{aligned}
R(\lambda, A^\Phi) &= R(\lambda, A) + L_\lambda \cdot \big(Id_{\partial X} - \Phi L_\lambda\big)^{-1} \cdot \Phi R(\lambda, A) \\
&= \big(Id_Z - L_\lambda \Phi\big)^{-1} \cdot R(\lambda, A).
\end{aligned}$$

*Remark* 3.14. Note that by definition $LL_\lambda = Id_{\partial X}$. Hence, the condition $1 \in \rho(\Delta_{\partial X})$ appearing in the previous result is equivalent to the fact that $\Psi L_\lambda$ is invertible where $\Psi := L - \Phi$ determines the domain $D(A^\Phi) = \ker(L - \Phi) = \ker(\Psi)$.

Corollaries 2.8 and 2.9 read in the context of the generic example as follows.

**Corollary 3.15.** *Let $\lambda_0 \in \partial\sigma(A)$.*

(a) *If the map $\rho(A) \ni \lambda \mapsto L_\lambda \Phi \in \mathcal{L}(Z)$ has a continuous extension in $\lambda_0$, then $\lambda_0 \in \sigma(A^\Phi)$.*

(b) *If for some $0 \neq f \in \ker(\lambda_0 - A)$ the limit*
$$\lim_{\rho(A) \ni \lambda \to \lambda_0} L_\lambda \Phi f =: g$$
*converges in $X$ such that $g \in Z$ and $\Phi g = 0$, then $\lambda_0 \in \sigma_p(A^\Phi)$.*

We mention that the operator $A^\Phi$ for bounded $\Phi \in \mathcal{L}(X, \partial X)$ was already studied by Greiner in [8]. In case $\dim(\partial X) < +\infty$, [8, Prop. 3.1] characterizes the spectral values of $A^\Phi$ lying in the component of $\rho(A)$ which is unbounded to the right as the zeros of the analytic function $F(\lambda) := \det(Id_{\partial X} - \Delta_{\partial X}(\lambda))$, cf. Corollary 3.13. Moreover, in [13] Nagel studied domain perturbations for operator matrices and arrives in [13, Thm. 2.7] at conditions similar to (3.12) assuming compactness of operators corresponding to our $\Delta_Z(\lambda) \in \mathcal{L}(Z)$.

In Section 4 we consider a series of concrete applications, most of which fit into the setting of the generic example above.

## 4. Applications

In this section we will apply our abstract results from Sections 2–3 to

(i) the first derivative with general boundary conditions,
(ii) the second derivative with general boundary conditions,
(iii) a second derivative with Nonlocal Neumann boundary conditions,
(iv) a second order differential operator with point delay at the boundary, and
(v) the Laplacian with dynamical boundary conditions.

Moreover, we use them to investigate the spectral theory of

(vi) delay equations,
(vii) complete second order Cauchy problems.

We note that in the examples (i)–(ii) a direct computation of the spectrum is also possible. Nevertheless, these examples illustrate in a simple context our results.

Finally, we mention that our approach can also be used for spectral investigations of flows in networks (cf. [11, Prop. 3.3]) and in various other situations like [10, Sect. II] and [14, Sect. 3].

### 4.1. The First Derivative with General B.Cs.
The aim of this example is to illustrate our results in a very simple but typical context. Let $X = \mathrm{L}^p[0,1]$ and $\Psi \in (\mathrm{W}^{1,p}[0,1])'$ for some $1 \le p < \infty$. We show how the operator

$$(4.1) \qquad G = \tfrac{d}{ds} \quad \text{with domain } D(G) := \left\{ f \in \mathrm{W}^{1,p}[0,1] : \Psi f = 0 \right\}$$

fits into the framework of our generic example from the previous section. In fact, it suffices to choose the maximal operator $A_m := \tfrac{d}{ds}$ with domain $D(A_m) = \mathrm{W}^{1,p}[0,1]$, the boundary space $\partial X := \mathbb{C}$ and $L = \delta_0 : D(A_m) \to \partial X$ where $\delta_0 f := f(0)$. This yields

$$A = \tfrac{d}{ds} \quad \text{with domain } D(A) := \left\{ f \in \mathrm{W}^{1,p}[0,1] : f(0) = 0 \right\} = \ker(L).$$

Moreover, let $Z := [D(A_m)] = \mathrm{W}^{1,p}[0,1]$, then $\Phi := L - \Psi \in \mathcal{L}(Z, \partial X)$ and by definition it follows that $G = A^\Phi$.

Next one easily verifies that $\rho(A) = \mathbb{C}$ and that for $\lambda \in \mathbb{C}$ the Dirichlet operators $L_\lambda \in \mathcal{L}(\mathbb{C}, X)$ are given by

$$L_\lambda z = z\, e^{\lambda \bullet}, \quad z \in \mathbb{C}.$$

Hence, Corollary 3.13 implies the following.

**Corollary 4.1.** *The spectrum of $G$ in* (4.1) *is characterized by*

$$\lambda \in \sigma(G) \iff \lambda \in \sigma_p(G) \iff \Psi\bigl(e^{\lambda \bullet}\bigr) = 0.$$

For example $\Psi = \delta_0 - \delta_1$ implies $\lambda \in \sigma(G) \iff e^\lambda = 1$, i.e., $\sigma(G) = \sigma_p(G) = 2\pi i \mathbb{Z}$.

*Remark* 4.2. We note that the choice of the unperturbed operator $A \subset A_m$ with domain $D(A) = \ker(L)$ in the example above (as well as in the following ones) is rather arbitrary. As already mentioned, due to the freedom of the perturbation $\Phi \in (\mathrm{W}^{1,p}[0,1])'$ it is convenient to choose $A$ having small spectrum to obtain the least possible points $\lambda \in \sigma(A)$ which have to be investigated separately using, e.g., Corollary 3.15. This fact is quite different from perturbation results for generators where in most cases the perturbation $\Phi$ has to be "small" in order that the generator property of $A$ is inherited to $A^\Phi$.

### 4.2. The Second Derivative with General B.Cs.
On the state space $X = \mathrm{C}[0,1]$ we consider for some $\psi_1, \psi_2 \in (\mathrm{C}^2[0,1])'$ the second derivative

$$(4.2) \qquad G = \tfrac{d^2}{ds^2} \quad \text{with domain } D(G) := \left\{ f \in \mathrm{C}^2[0,1] : \psi_1(f) = 0 = \psi_2(f) \right\}.$$

To compute $\sigma(G)$ we consider the maximal operator $A_m := \tfrac{d^2}{ds^2}$ with domain $D(A_m) := \mathrm{C}^2[0,1]$, the boundary space $\partial X = \mathbb{C}^2$ and $L := \binom{\delta_0}{\delta_0'} : D(A_m) \to \partial X$ where $\delta_0' f := f'(0)$. This gives the second derivative

$$A = \tfrac{d^2}{ds^2} \quad \text{with domain } D(A) := \left\{ f \in \mathrm{C}^2[0,1] : Lf = 0 \right\} = \ker(L).$$

Moreover, for $Z := [D(A_m)] = \mathrm{C}^2[0,1]$ we have $\Phi = \binom{\varphi_1}{\varphi_2} := L - \binom{\psi_1}{\psi_2} \in \mathcal{L}(Z, \partial X)$ and by definition it follows $G = A^\Phi$.

Since by Arzela–Ascoli's theorem the embedding $[D(A)] \hookrightarrow X$ is compact, the operator $A$ has compact resolvent, which implies $\sigma(A) = \sigma_p(A)$. Now a simple computations shows that

$$\sigma_p(A) = \emptyset.$$

Next, by solving for $\lambda \in \mathbb{C}$ and $x := \binom{z_1}{z_2} \in \partial X = \mathbb{C}^2$ the Dirichlet problem

$$\begin{cases} (\lambda - A_m)f = 0, \\ Lf = x \end{cases} \iff \begin{cases} (\lambda - \tfrac{d^2}{ds^2})f = 0, \\ f(0) = z_1, \ f'(0) = z_2 \end{cases}$$

we obtain the Dirichlet operators $L_\lambda \in \mathcal{L}(\partial X, X) = \mathcal{L}(\mathbb{C}^2, \mathrm{C}[0,1])$ given by

$$\left(L_\lambda \binom{z_1}{z_2}\right)(s) = \begin{cases} z_1 \cdot \cosh(\sqrt{\lambda}s) + \tfrac{z_2 \cdot \sinh(\sqrt{\lambda}s)}{\sqrt{\lambda}} & \text{if } \lambda \neq 0, \\ z_1 + z_2 \cdot s & \text{if } \lambda = 0, \end{cases}$$

for $\binom{z_1}{z_1} \in \partial X$ and $s \in [0,1]$. Now Corollary 3.13 applied to this situation gives the following.

**Corollary 4.3.** *For $G$ given by* (4.2) *we have $\lambda \in \sigma(G) = \sigma_p(G)$ if and only if*

$$(4.3) \qquad \begin{cases} \det \begin{pmatrix} \psi_1(\cosh(\sqrt{\lambda}\bullet)), & \psi_1(\sinh(\sqrt{\lambda}\bullet)) \\ \psi_2(\cosh(\sqrt{\lambda}\bullet)), & \psi_2(\sinh(\sqrt{\lambda}\bullet)) \end{pmatrix} = 0 & \text{for } \lambda \neq 0, \\ \det \begin{pmatrix} \psi_1(\mathbb{1}) & \psi_1(\mathbb{s}) \\ \psi_2(\mathbb{1}) & \psi_2(\mathbb{s}) \end{pmatrix} = 0 & \text{for } \lambda = 0, \end{cases}$$

*where $\mathbb{1}(s) = 1$ and $\mathbb{s}(s) = s$ for all $s \in [0,1]$.*

For particular choices of the boundary functionals $\psi_1, \psi_2$ the characteristic equation (4.3) might simplify considerably. For example, if we consider the second derivative $G = \tfrac{d^2}{ds^2}$ with Wentzell-type boundary conditions $f''(j) = f'(j)$, $j = 0, 1$, we obtain the following.

**Corollary 4.4.** *For $\Psi = \binom{\delta_0'' - \delta_0'}{\delta_1'' - \delta_1'}$ we obtain $\sigma(G) = \sigma_p(G) = \{-n^2 \cdot \pi^2 : n \in \mathbb{N}_0\} \cup \{1\}$.*

*Proof.* By (4.3) we have $\lambda \in \sigma(G)$ if and only if

$$\lambda \cdot (\lambda - 1) \cdot \sinh(\sqrt{\lambda}) = 0. \qquad \square$$

### 4.3. The Second Derivative with Nonlocal Neumann B.Cs.
This example deals with an operator associated to a heat equation with distributed unbounded and delayed feedback at the boundary, cf. [9, 2].

For $1 \leq p < +\infty$ let $X := \mathrm{L}^p[0,1]$ and $Y := \mathrm{L}^p([-1,0], X)$, which by [4, Thm. A.6] is isometrically isomorphic to $\mathrm{L}^p([-1,0] \times [0,1])$. For this reason in the sequel we will use the notation $v(r,s) :=$ $(v(r))(s)$ for $v \in Y$ and $r \in [-1,0]$, $s \in [0,1]$. Then on the product space $\mathcal{X} := X \times Y$ we consider the operator matrix

$$\mathcal{G} := \begin{pmatrix} \frac{d^2}{ds^2} & 0 \\ 0 & \frac{d}{dr} \end{pmatrix},$$

$$D(\mathcal{G}) := \Big\{ \begin{pmatrix} f \\ v \end{pmatrix} \in \mathrm{W}^{2,p}[0,1] \times \mathrm{W}^{1,p}\big([-1,0], X\big) : v(0) = f,\ f(1) = 0,$$
$$f'(0) = \int_0^1 \int_{-1}^0 v(r,s)\, d\nu(r)\, ds \Big\}$$

where $\nu : [-1,0] \to \mathbb{R}$ is a function of bounded variation. This operator appears in [9, Expl. 5.2] (for $p = 2$) and in [2, Sect. 3.3], where it is shown that it generates a $C_0$-semigroup.

Here we are interested in characterizing the spectrum of $\mathcal{G}$. In order to represent $\mathcal{G}$ as $\mathcal{A}^\Phi$ we first introduce the following operators and spaces. Consider

- $A_m := \frac{d^2}{ds^2}$ with domain $D(A_m) = \{f \in \mathrm{W}^{2,p}[0,1] : f(1) = 0\}$ on $X$,
- $L := \delta'_1 : D(A_m) \to \partial X := \mathbb{C}$, i.e., $Lf = f'(1)$,
- $D_m := \frac{d}{dr}$ with domain $D(D_m) = \mathrm{W}^{1,p}([-1,0], X)$ on $Y$,
- $K := \delta_0 : D(D_m) \to \partial Y := X$, i.e., $Kv = v(0)$,
- $A = A_m|_{\ker L}$, $D := D_m|_{\ker K}$.

Next we define the maximal operator matrix

$$\mathcal{A}_m := \begin{pmatrix} A_m & 0 \\ 0 & D_m \end{pmatrix}, \quad D(\mathcal{A}_m) := D(A_m) \times D(D_m).$$

Moreover, let $\partial \mathcal{X} := \partial X \times \partial Y$,

$$\mathcal{L} := \begin{pmatrix} L & 0 \\ 0 & K \end{pmatrix} : D(\mathcal{A}_m) \to \partial \mathcal{X}$$

and $\mathcal{A} \subset \mathcal{A}_m$ with domain $D(\mathcal{A}) := \ker(\mathcal{L}) = D(A) \times D(D)$. Finally, we take $\mathcal{Z} := X \times [D(D_m)]$ and define

$$\Phi := \begin{pmatrix} \delta'_1 - \delta'_0 & \varphi \\ Id_X & 0 \end{pmatrix} \in \mathcal{L}(\mathcal{Z}, \partial \mathcal{X}) \quad \text{where} \quad \varphi(v) := \int_0^1 \int_{-1}^0 v(r,s)\, d\nu(r)\, ds.$$

Then by definition we obtain $\mathcal{G} = \mathcal{A}^\Phi$. In order to characterize the spectrum of $\mathcal{A}^\Phi$ we first note that $\sigma(A) = \sigma(D) = \emptyset$. Moreover, for $\lambda \in \mathbb{C}$ the Dirichlet operators $L_\lambda \in \mathcal{L}(\partial X, X)$ and $K_\lambda \in \mathcal{L}(\partial Y, Y)$ relatively to the pairs $A_m, L$ and $D_m, K$ are given by

$$(L_\lambda z)(s) = \begin{cases} z \cdot \frac{\sinh(\sqrt{\lambda}(s-1))}{\sqrt{\lambda}} & \text{if } \lambda \neq 0, \\ z \cdot (s-1) & \text{if } \lambda = 0, \end{cases} \qquad z \in \partial X,\ s \in [0,1],$$

$$(K_\lambda f)(r) = e^{\lambda r} \cdot f, \qquad\qquad\qquad\qquad\qquad\qquad\qquad f \in \partial Y,\ r \in [-1,0].$$

Thus, for $\lambda \in \mathbb{C}$ we obtain the Dirichlet operator for the pair $\mathcal{A}_m, \mathcal{L}$ as

$$\mathcal{L}_\lambda := \begin{pmatrix} L_\lambda & 0 \\ 0 & K_\lambda \end{pmatrix} \in \mathcal{L}(\partial \mathcal{X}, \mathcal{X}).$$

Now we are in the position to apply Corollary 3.13 and obtain the following characterization of the spectral values of $\mathcal{G} = \mathcal{A}^\Phi$.

**Corollary 4.5.** *Let $\lambda \in \mathbb{C}$. Then for $l_\lambda := L_\lambda 1$ we have*

$$\lambda \in \sigma(\mathcal{A}^\Phi) = \sigma_p(\mathcal{A}^\Phi) \quad \Longleftrightarrow \quad \int_0^1 \int_{-1}^0 e^{\lambda r} \cdot l_\lambda(s)\, d\nu(r)\, ds = \cosh(\sqrt{\lambda}).$$

*In particular, if $\nu = \delta_{-1}$ then*

$$\lambda \in \sigma_p(\mathcal{A}^\Phi) \quad \Longleftrightarrow \quad \big(\lambda \cdot e^\lambda + 1\big) \cdot \cosh(\sqrt{\lambda}) = 1.$$

*Proof.* For $\lambda \in \mathbb{C}$ we have
$$\Phi\mathcal{L}_\lambda = \begin{pmatrix} 1 - \cosh(\sqrt{\lambda}) & \varphi K_\lambda \\ L_\lambda & 0 \end{pmatrix} \in \mathcal{L}(\partial\mathcal{X}) = \mathcal{L}\big(\mathbb{C} \times \mathrm{L}^p[0,1]\big).$$
By Corollary 3.13 and Lemma A.10 this implies that
$$\begin{aligned}
\lambda \in \sigma(\mathcal{A}^\Phi) &\implies 1 \in \sigma(\Phi\mathcal{L}_\lambda) &\iff& \quad 0 \in \sigma\big(\cosh(\sqrt{\lambda}) - \varphi K_\lambda L_\lambda\big) \\
&\iff 0 \in \sigma_p\big(\cosh(\sqrt{\lambda}) - \varphi K_\lambda L_\lambda\big) &\iff& \quad 1 \in \sigma_p(\Phi\mathcal{L}_\lambda) \\
&\iff \lambda \in \sigma_p(\mathcal{A}^\Phi) &\implies& \quad \lambda \in \sigma(\mathcal{A}^\Phi),
\end{aligned}$$
where we used that $\dim(\partial X)$ is finite (in fact one) dimensional. The assertions then follows by computing $\varphi K_\lambda L_\lambda : \mathbb{C} \to \mathbb{C}$. □

4.4. **A Second Order Differential Operator with Point Delay at the Boundary.** This example treats an operator associated to a one-dimensional reaction-diffusion equation modeling a delayed chemical reaction.

Let $X := \mathrm{L}^p[0,1]$ and $Y := \mathrm{L}^p([-1,0], \partial Y)$ for $\partial Y := \mathrm{W}^{1,p}[0,1]$. Then on the product space $\mathcal{X} := X \times Y$ we consider for some fixed $c, k \in \mathbb{C}$ the operator matrix
$$\mathcal{G} := \begin{pmatrix} \frac{d^2}{ds^2} - 2c \cdot \frac{d}{ds} + k \cdot Id_X & 0 \\ 0 & \frac{d}{dr} \end{pmatrix},$$
$$D(\mathcal{G}) := \left\{ \begin{pmatrix} f \\ v \end{pmatrix} \in \mathrm{W}^{2,p}[0,1] \times \mathrm{W}^{1,p}([-1,0], \partial Y) : \begin{matrix} \binom{f'(0)}{f'(1)} = \binom{f(0)}{0} - \binom{v(-1,1)}{0} \\ v(0) = f \end{matrix} \right\}.$$

Using [2, Cors. 3.6 & 3.7] one can show that $\mathcal{G}$ generates a $C_0$-semigroup on $\mathcal{X}$. In order to compute $\sigma(\mathcal{G})$ we introduce the following operators and spaces.
- $A_m := \frac{d^2}{ds^2} - 2c \cdot \frac{d}{ds} + k \cdot Id_X$ with domain $D(A_m) = \{f \in \mathrm{W}^{2,p}[0,1] : f'(1) = 0\}$ on $X$,
- $L := \delta_1 : D(A_m) \to \partial X := \mathbb{C}$, i.e., $Lf = f(1)$,
- $D_m := \frac{d}{dr}$ with domain $D(D_m) = \mathrm{W}^{1,p}([-1,0], \partial Y)$ on $Y$,
- $K := \delta_0 : D(D_m) \to \partial Y$, i.e., $Kv = v(0)$,
- $A = A_m|_{\ker L}$, $D := D_m|_{\ker K}$.

Next we define the maximal operator matrix
$$\mathcal{A}_m := \begin{pmatrix} A_m & 0 \\ 0 & D_m \end{pmatrix}, \quad D(\mathcal{A}_m) := D(A_m) \times D(D_m).$$
Moreover, let $\partial\mathcal{X} := \partial X \times \partial Y$,
$$\mathcal{L} := \begin{pmatrix} L & 0 \\ 0 & K \end{pmatrix} : D(\mathcal{A}_m) \to \partial\mathcal{X}$$
and $\mathcal{A} \subset \mathcal{A}_m$ with domain $D(\mathcal{A}) := \ker(\mathcal{L}) = D(A) \times D(D)$. Finally, we define the spaces $Z_X := \mathrm{W}^{1,p}[0,1]$, $Z_Y := \mathrm{W}^{1,p}([-1,0], \partial Y)$, $\mathcal{Z} := Z_X \times Z_Y$ and consider
$$\Phi := \begin{pmatrix} \varphi & \psi \\ Id_{Z_X} & 0 \end{pmatrix} \in \mathcal{L}(\mathcal{Z}, \partial\mathcal{X}), \quad \text{where}$$
$$\varphi := \delta_0 + \delta_1 - \delta_0' \in \mathcal{L}(Z_X, \partial X), \quad \text{i.e., } \varphi(f) = f(0) + f(1) - f'(0),$$
$$\psi := -\delta_{-1} \otimes \delta_1 \in \mathcal{L}(Z_Y, \partial X), \quad \text{i.e., } \psi(v) = -v(-1,1).$$

Then by definition we obtain $\mathcal{G} = \mathcal{A}^\Phi$. In order to characterize the spectrum of $\mathcal{A}^\Phi$ we first note that $\sigma(A) = \sigma(D) = \emptyset$. Moreover, the Dirichlet operators $L_\lambda \in \mathcal{L}(\partial X, X)$ and $K_\lambda \in \mathcal{L}(\partial Y, Y)$ relatively to the pairs $A_m, L$ and $D_m, K$ are explicitly given by
$$(L_\lambda z)(s) = \begin{cases} e^{c(s-1)} \cdot \left( \cosh\big((s-1) \cdot \sqrt{\lambda + c^2 - k}\big) - \frac{c \sinh\big((s-1)\cdot\sqrt{\lambda+c^2-k}\big)}{\sqrt{\lambda+c^2-k}} \right) \cdot z & \text{if } \lambda \neq k - c^2, \\ e^{c(s-1)} \cdot (1 + c - cs) \cdot z & \text{if } \lambda = k - c^2, \end{cases}$$
$$(K_\lambda f)(r) = e^{\lambda r} \cdot f,$$

where $z \in \partial X = \mathbb{C}$, $s \in [0,1]$, $f \in \partial Y = W^{1,p}[0,1]$ and $r \in [-1, 0]$. Thus, for $\lambda \in \mathbb{C}$ we obtain the Dirichlet operator for the pair $\mathcal{A}_m$, $\mathcal{L}$ as

$$\mathcal{L}_\lambda := \begin{pmatrix} L_\lambda & 0 \\ 0 & K_\lambda \end{pmatrix} \in \mathcal{L}(\partial \mathcal{X}, \mathcal{X}).$$

Now by Corollary 3.13 we obtain the following characterization of the spectral values of $\mathcal{G} = \mathcal{A}^\Phi$.

**Corollary 4.6.** *Let $\lambda \in \mathbb{C}$. Then we have $\lambda \in \sigma(\mathcal{G}) = \sigma_p(\mathcal{G})$ if and only if*

(4.4) $$e^{-\lambda} - l_\lambda(0) + l'_\lambda(0) = 0,$$

*where $l_\lambda(s) := (L_\lambda 1)(s)$, $s \in [0,1]$.*

*Proof.* For arbitrary $\lambda \in \mathbb{C}$ we have

$$Id_{\partial \mathcal{X}} - \Phi \mathcal{L}_\lambda = \begin{pmatrix} l'_\lambda(0) - l_\lambda(0) & e^{-\lambda} \cdot \delta_1 \\ -L_\lambda & Id_{\partial Y} \end{pmatrix} \in \mathcal{L}(\partial \mathcal{X}).$$

Using Schur complements from Appendix A.4, this matrix is not invertible if and only if (4.4) holds. The assertion then follows from Corollary 3.13. □

**4.5. The Laplacian with Dynamical Boundary Conditions.** Let $X := L^2(\Omega)$ and $Y := L^2(\partial \Omega)$ for some open, bounded domain $\Omega \subset \mathbb{R}^n$ with smooth boundary $\partial \Omega$. We consider the Laplace operator $\Delta_\Omega$ on $X$ and the Laplace-Beltrami operator $\Delta_{\partial \Omega}$ on $Y$ with domains

$$D(\Delta_\Omega) := \left\{ f \in H^{\frac{3}{2}}(\Omega) \cap H^2_{loc}(\Omega) : \Delta f \in L^2(\Omega) \right\},$$

$$D(\Delta_{\partial \Omega}) := \left\{ f \in L^2(\partial \Omega) : \Delta f \in L^2(\partial \Omega) \right\}.$$

Then as in [5, Sects. 3 & 6] we define on the space $\mathcal{X} := X \times Y$ the operator matrix

(4.5) $$\mathcal{G} := \begin{pmatrix} \Delta_\Omega & 0 \\ C & \Delta_{\partial \Omega} \end{pmatrix},$$
$$D(\mathcal{G}) := \left\{ \begin{pmatrix} f \\ g \end{pmatrix} \in D(\Delta_\Omega) \times D(\Delta_{\partial \Omega}) : \frac{\partial f}{\partial \nu}\big|_{\partial \Omega} = g \right\}$$

for some operator $C \in \mathcal{L}(H^1(\Omega), L^2(\partial \Omega))$, e.g., $Cg = \int_\Omega k(s, \bullet) \nabla g(s)\, ds \in L^2(\partial \Omega)$ for a $L^2$-function $k : \Omega \times \partial \Omega \to \mathbb{C}^n$.

In order to embed $\mathcal{G}$ in our setting we could either write it as $\mathcal{G} = \mathcal{A}^\Phi_\mathcal{P}$ like we did in the previous examples or represent it as $\mathcal{G} = \mathcal{A}_{\mathcal{B}\mathcal{C}}$. Here we use the second alternative and introduce to this end the following operators and spaces.

- The operator $L : D(\Delta_\Omega) \to L^2(\partial \Omega)$, $Lf := \frac{\partial f}{\partial \nu}\big|_{\partial \Omega}$,
- the Neumann Laplacian $A := \Delta_N \subset \Delta_\Omega$ with domain $D(A) := \ker L$ which satisfies $\sigma(A) \subset \{\zeta \in \mathbb{C} : \text{Re}(\zeta) \le 0\}$,
- $Z := H^1(\Omega)$ and $Z_{-1} := (\eta_0 - A_{-1})Z$ for some $\eta_0 > 0$, cf. Section A.2,
- $A^Z := A_{-1}|_{Z_{-1}} : Z \subset Z_{-1} \to Z_{-1}$,
- $\mathcal{Z} := H^1(\Omega) \times D(\Delta_{\partial \Omega})$, $\mathcal{Z}_{-1} : Z_{-1} \times Y$ and $\mathcal{U} := Y \times Y$.

Then we define the operator matrices

$$\mathcal{A}^\mathcal{Z} := \begin{pmatrix} A^Z & 0 \\ 0 & \Delta_{\partial \Omega} \end{pmatrix} : \mathcal{Z} \subset \mathcal{Z}_{-1} \to \mathcal{Z}_{-1}$$

and

$$\mathcal{B} := \begin{pmatrix} L_A & 0 \\ 0 & Id_Y \end{pmatrix} \in \mathcal{L}(\mathcal{U}, \mathcal{X}^\mathcal{A}_{-1}), \qquad \mathcal{C} := \begin{pmatrix} 0 & Id_Y \\ C & 0 \end{pmatrix} \in \mathcal{L}(\mathcal{Z}, \mathcal{U}).$$

Here $L_A := (\eta_0 - A^Z)L_{\eta_0} \in \mathcal{L}(Y, Z_{-1})$ where the Dirichlet operator $L_{\eta_0} \in \mathcal{L}(Y, Z)$ exists by [5, p.298–299]. Then one easily verifies that $\mathcal{G} = \mathcal{A}_{\mathcal{B}\mathcal{C}}$. By applying Theorem 2.3 we obtain the following characterization of the spectral values of $\mathcal{G}$.

**Corollary 4.7.** *For $\mathcal{G}$ given by (4.5) and $\lambda \in \rho(\mathcal{A}) = \rho(\Delta_N) \cap \rho(\Delta_{\partial \Omega})$ we have*

$$\lambda \in \sigma(\mathcal{G}) \iff \lambda \in \sigma_p(\mathcal{G}) \iff \lambda \in \sigma_p(\Delta_{\partial \Omega} + CL_\lambda).$$

*Proof.* Let $\lambda \in \rho(\mathcal{A}) = \rho(A) \cap \rho(\Delta_{\partial\Omega})$. Then by (3.4) we have $R(\lambda, \mathcal{A}^Z)L_A = L_\lambda$ which implies

$$\Delta_{\mathcal{U}}(\lambda) = \mathcal{C}R(\lambda, \mathcal{A}^Z)\mathcal{B} = \begin{pmatrix} 0 & R(\lambda, \Delta_{\partial\Omega}) \\ CL_\lambda & 0 \end{pmatrix} \in \mathcal{L}(Y \times Y).$$

By Theorem 2.3 and Lemma A.10 this gives

$$
\begin{aligned}
\lambda \in \sigma(\mathcal{G}) &\implies 1 \in \sigma(\Delta_{\mathcal{U}}(\lambda)) \\
&\iff 1 \in \sigma\bigl(CL_\lambda R(\lambda, \Delta_{\partial\Omega})\bigr) = \sigma_p\bigl(CL_\lambda R(\lambda, \Delta_{\partial\Omega})\bigr) \\
&\iff \lambda \in \sigma_p(\Delta_{\mathcal{U}}) \iff \lambda \in \sigma_p(\mathcal{G}), \\
&\implies \lambda \in \sigma(\mathcal{G}),
\end{aligned}
\tag{4.6}
$$

where the equality in (4.6) holds since the operator $\Delta_{\partial\Omega}$ has compact resolvent, yielding the compactness of $R(\lambda, \Delta_{\partial\Omega})CL_\lambda$. This proves that all the above conditions are equivalent. Since

$$\lambda - \Delta_{\partial\Omega} - CL_\lambda = \bigl(Id_Y - CL_\lambda R(\lambda, \Delta_{\partial\Omega})\bigr) \cdot \bigl(\lambda - \Delta_{\partial\Omega}\bigr)$$

where $\lambda - \Delta_{\partial\Omega} : D(\Delta_{\partial\Omega}) \to Y$ is bijective, we finally conclude that $\lambda \in \sigma_p(\Delta_{\partial\Omega} + CL_\lambda)$ iff $1 \in \sigma(CL_\lambda R(\lambda, \Delta_{\partial\Omega}))$ as claimed. $\square$

**4.6. Spectral Theory for Delay Equations.** In this section we will apply our results to delay equations of the form

(DE) $\qquad \begin{cases} \dot{u}(t) = Au(t) + \varphi \, u_t, & t \geq 0, \\ u(0) = x, \\ u(t) = f(t), & t \in [-1, 0]. \end{cases}$

Here $u : [-1, +\infty) \to X$ is a Banach space valued function, $u_t(\bullet) := u(\bullet + t)$, $A, \varphi$ are linear operators and $x, f$ are given initial values. There exist several approaches to (DE) depending on the underlying function space. Here we will consider the treatment on the spaces of $p$-integrable and continuous functions.

*4.6.a. Spectral Theory for the Reduction Matrix in $\mathrm{L}^p([-1,0], Y)$.* The $\mathrm{L}^p$-approach to delay equations (DE) yields to the operator matrix

(4.7) $\quad \mathcal{G} := \begin{pmatrix} A & \varphi \\ 0 & \frac{d}{ds} \end{pmatrix} \quad$ with domain $\quad D(\mathcal{G}) := \left\{ \begin{pmatrix} x \\ f \end{pmatrix} \in D(A) \times \mathrm{W}^{1,p}([-1,0], Y) : f(0) = x \right\}$

on $\mathcal{X} := X \times \mathrm{L}^p([-1,0], Y)$ for some $1 \leq p < +\infty$, cf. [4, Chap. 3]. Here

- $A : D(A) \subseteq X \to X$ is a linear operator on a Banach space $X$ satisfying $\rho(A) \neq \emptyset$,
- $Y$ is a Banach space such that $[D(A)] \hookrightarrow Y \hookrightarrow X$,
- $\varphi \in \mathcal{L}(\mathrm{W}^{1,p}([-1,0], Y), X)$.

In order to represent $\mathcal{G}$ as in our generic example we introduce the following operators and spaces. By $\mathrm{diag}(\ldots)$ we denote a diagonal matrix with the given entries.

- $\mathcal{A}_m := \mathrm{diag}(A, \frac{d}{ds}) : D(\mathcal{A}_m) \subset \mathcal{X} \to \mathcal{X}$ with domain $D(\mathcal{A}_m) := D(A) \times \mathrm{W}^{1,p}([-1,0], Y)$,
- $\mathcal{Z} := [D(\mathcal{A}_m)]$, $\partial \mathcal{X} := Y$,
- $\mathcal{L} = (0, \delta_0) \in \mathcal{L}(\mathcal{Z}, \partial\mathcal{X})$, i.e., $\mathcal{L}(x, f)^\top = f(0)$,
- $\mathcal{A} := \mathcal{A}_m|_{\ker(\mathcal{L})}$,
- $\mathcal{P} := \begin{pmatrix} 0 & \varphi \\ 0 & 0 \end{pmatrix} \in \mathcal{L}(\mathcal{Z}, \mathcal{X})$,
- $\Phi = (Id_Y, 0) \in \mathcal{L}(\mathcal{Z}, \partial\mathcal{X})$ i.e., $\Phi(x, f)^\top = x$.

Then $\rho(\mathcal{A}) = \rho(A) \neq \emptyset$, $\mathcal{L}$ is surjective and $\mathcal{A}_m$ is closed, hence the Assumptions 3.5 are satisfied. Moreover, one easily verifies that

$$\mathcal{G} = \mathcal{A}_\mathcal{P}^\Phi.$$

From the results of Subsection 3.4 we now obtain the following spectral characterizations, cf. [4, Prop. 3.19, Lem. 3.20].

**Corollary 4.8.** *For all $\lambda \in \mathbb{C}$ and $1 \leq p < +\infty$ we have*
$$\lambda \in \sigma(\mathcal{G}) \iff \lambda \in \sigma(A + \varphi L_\lambda),$$
$$\lambda \in \sigma_p(\mathcal{G}) \iff \lambda \in \sigma_p(A + \varphi L_\lambda),$$
*where $L_\lambda \in \mathcal{L}(X, W^{1,p}([-1,0], Y))$ is given by $(L_\lambda x)(s) := e^{\lambda s} x$ for $x \in X$ and $s \in [-1, 0]$. Moreover, if $\rho(\mathcal{G}) \neq \emptyset$ then*
$$\lambda \in \sigma_*(\mathcal{G}) \iff \lambda \in \sigma_*(A + \varphi L_\lambda)$$
*for all $* \in \{a, r, c, ess\}$.*

*Proof.* We apply Corollary 3.10. To this end fix some $\lambda_0 \in \rho(A)$. Then for $\lambda \in \mathbb{C}$ we obtain
$$\Delta_{\mathcal{Z}}(\lambda_0, \lambda) = \begin{pmatrix} (\lambda_0 - \lambda) R(\lambda_0, A) & R(\lambda_0, A) \varphi \\ L_{\lambda_0} & (\lambda_0 - \lambda) R(\lambda_0, D) \end{pmatrix} \in \mathcal{L}(\mathcal{Z})$$
for $D := \frac{d}{ds}$ with domain $D(D) := W_0^{1,p}([-1,0], Y) = \ker(\delta_0)$ which satisfies $\rho(D) = \mathbb{C}$. Since $\mathcal{A}_m + \mathcal{P}$ is closed and $\mathcal{L} - \Phi$ is surjective, by (3.9) we conclude that

$$\begin{aligned}
\lambda \in \sigma(\mathcal{G}) &\iff 1 \in \rho(\Delta_{\mathcal{Z}}(\lambda_0, \lambda)) \\
&\iff \begin{pmatrix} Id - (\lambda_0 - \lambda) R(\lambda_0, A) & -R(\lambda_0, A)\varphi \\ -L_{\lambda_0} & (\lambda - D) R(\lambda_0, D) \end{pmatrix} \quad \text{is invertible in } \mathcal{L}(\mathcal{Z}) \\
&\iff Id - R(\lambda_0, A) \cdot (\lambda_0 - \lambda + \varphi \cdot (\lambda_0 - D) R(\lambda, D) \cdot L_{\lambda_0}) \quad \text{is invertible in } \mathcal{L}([D(A)]) \\
&\iff Id - (\lambda_0 - \lambda + \varphi \cdot L_\lambda) \cdot R(\lambda_0, A) \quad \text{is invertible in } \mathcal{L}(X) \\
&\iff \lambda - A - \varphi \cdot L_\lambda : D(A) \subseteq X \to X \quad \text{is invertible in } \mathcal{L}(X),
\end{aligned}$$

where we used the resolvent equation, the Schur complement in $[D(A)]$, (3.4) and Corollary A.12. The assertions concerning the subdivision of the spectrum then follows from Corollary 3.9.(d). □

*Remark* 4.9. By [4, Thm. 3.12] the delay equation (DE) is well posed if and only if $\mathcal{G}$ generates a $C_0$-semigroup on $\mathcal{X}$. In this case $\rho(\mathcal{G}) \neq \emptyset$, hence the previous result gives a complete description of the spectrum and its finer subdivisions.

4.6.b. *Spectral Theory for the First Derivative in* $C([-1,0], X)$. If one treats the delay equation (DE) within a framework of continuous functions then its initial values $x, f$ always satisfy $f(0) = x$, i.e. this condition is superfluous. For this reason the reduction to an abstract Cauchy problem does not yield an operator matrix as in the previous subsection. Instead, one obtains the operator

(4.8) $\quad \mathcal{G} := \frac{d}{ds} \quad \text{with domain} \quad D(\mathcal{G}) := \{f \in C^1([-1,0], X) : f(0) \in D(A), \ f'(0) = Af(0) + \varphi f\}$

on $\mathcal{X} := C([-1,0], X)$ equipped with the sup-norm $\|\bullet\|_\infty$, cf. [7, Sect. VI.6]. Here we assume that
- $A : D(A) \subseteq X \to X$ is a linear operator on a Banach space $X$ satisfying $\rho(A) \neq \emptyset$,
- $\varphi \in \mathcal{L}(\mathcal{Z}, X)$ for $\mathcal{Z} := C^1([-1,0], X)$.

In order to represent $\mathcal{G}$ as in our generic example we first observe that for $\lambda_0 \in \rho(A)$ we have $f \in D(\mathcal{G})$ if and only if

(4.9) $\qquad\qquad f(0) = R(\lambda_0, A)(\varphi f + \lambda_0 f(0) - f'(0)).$

Next we introduce the following operators and spaces.
- $\mathcal{A}_m := \frac{d}{ds} : D(\mathcal{A}_m) \subset \mathcal{X} \to \mathcal{X}$ with domain $D(\mathcal{A}_m) := C^1([-1,0], X)$,
- $\mathcal{L} = \delta_0 \in \mathcal{L}(\mathcal{Z}, \partial\mathcal{X})$, i.e., $Lf = f(0)$, for $\partial\mathcal{X} := X$,
- $\mathcal{A} := \mathcal{A}_m|_{\ker(\mathcal{L})}$,
- $\Phi := R(\lambda_0, A)(\varphi + \lambda_0 \delta_0 - \delta'_0) \in \mathcal{L}(\mathcal{Z}, \partial\mathcal{X})$ i.e., $\Phi f = R(\lambda_0, A)(\varphi f + \lambda_0 f(0) - f'(0))$.

Then $\rho(\mathcal{A}) = \mathbb{C}$, $\mathcal{L}$ is surjective and $\mathcal{A}_m$ is closed, hence the Assumptions 3.5 are satisfied. Moreover, using (4.9), one easily verifies that
$$\mathcal{G} = \mathcal{A}^\Phi.$$

Before applying the results from Subsection 3.4 to $\mathcal{G}$, we study the surjectivity of $\mathcal{L} - \Phi$. Here we define for $\mu \in \mathbb{C}$ the functions $\varepsilon_\mu \in C^1[-1,0]$ by $\varepsilon_\mu(s) := e^{\mu s}$.

**Lemma 4.10.** *If there exists $\lambda_0 \in \rho(A)$ and $\mu \in \mathbb{C}$, $\mu \neq \lambda_0$, such that*
$$\left\|\varphi\bigl((\varepsilon_{\lambda_0} - \varepsilon_\mu) \otimes Id\bigr)\right\|_{\mathcal{L}(X)} < |\lambda_0 - \mu|, \tag{4.10}$$
*then $\mathcal{L} - \Phi$ is surjective. This is the case if there exist $\mu_n \in \mathbb{C}$, $n \in \mathbb{N}$, such that $|\mu_n| \to +\infty$ and*
$$\lim_{n \to +\infty} \frac{\left\|\varphi(\varepsilon_{\mu_n} \otimes Id)\right\|_{\mathcal{L}(X)}}{|\mu_n|} < 1.$$
*In particular, if $\varphi \in \mathcal{L}(\mathcal{X}, X)$ then $\mathcal{L} - \Phi$ is always surjective.*

*Proof.* We show that (4.10) implies the inclusion (3.10), hence $\mathcal{L} - \Phi$ is surjective by Lemma 3.12. In fact, since $\ker(\lambda_0 - \mathcal{A}_m) = \varepsilon_{\lambda_0} \otimes X$ and $(\varepsilon_{\lambda_0} - \varepsilon_\mu) \otimes X \subset D(A)$, (3.10) follows if for each $x \in X$ there exists $y \in X$ such that
$$\begin{aligned}
& \Phi(\varepsilon_{\lambda_0} \otimes x) = \Phi\bigl((\varepsilon_{\lambda_0} - \varepsilon_\mu) \otimes y\bigr) \\
\iff \quad & \varphi(\varepsilon_{\lambda_0} \otimes x) = \varphi\bigl((\varepsilon_{\lambda_0} - \varepsilon_\mu) \otimes y\bigr) - (\lambda_0 - \mu) \cdot y \\
\iff \quad & \left(Id - \frac{\varphi\bigl((\varepsilon_{\lambda_0} - \varepsilon_\mu) \otimes Id\bigr)}{\lambda_0 - \mu}\right) \cdot y = \frac{\varphi(\varepsilon_{\lambda_0} \otimes x)}{\mu - \lambda_0} \\
\iff \quad & y = \left(Id - \frac{\varphi\bigl((\varepsilon_{\lambda_0} - \varepsilon_\mu) \otimes Id\bigr)}{\lambda_0 - \mu}\right)^{-1} \cdot \frac{\varphi(\varepsilon_{\lambda_0} \otimes x)}{\mu - \lambda_0},
\end{aligned}$$
where the invertibility of the operator in the last equivalence follows by assumption (4.10). The remaining two assertions follow easily by considering (4.10) for $\mu = \mu_n$ for sufficiently big $n \in \mathbb{N}$. □

Combining the previous result with Corollary 3.13 we immediately obtain the following spectral characterizations which significantly generalizes [7, Prop. VI.6.7].

**Corollary 4.11.** *For $\lambda \in \mathbb{C}$ define $L_\lambda \in \mathcal{L}(X, \mathrm{C}^1([-1,0], X))$ by $(L_\lambda x)(s) := e^{\lambda s} x$ for $x \in X$ and $s \in [-1, 0]$. Then for $\mathcal{G}$ defined in (4.8) we have*
$$\lambda \in \sigma_p(\mathcal{G}) \quad \iff \quad \lambda \in \sigma_p(A + \varphi L_\lambda).$$
*Moreover, if (4.10) is satisfied then*
$$\lambda \in \sigma(\mathcal{G}) \quad \iff \quad \lambda \in \sigma(A + \varphi L_\lambda).$$
*Finally, if in addition $\rho(\mathcal{G}) \neq \emptyset$ then*
$$\lambda \in \sigma_*(\mathcal{G}) \quad \iff \quad \lambda \in \sigma_*(A + \varphi L_\lambda)$$
*for all $* \in \{a, r, c, ess\}$.*

*Remark 4.12.* By [7, Cor. VI.6.3] the delay equation (DE) is well posed if and only if $\mathcal{G}$ generates a $C_0$-semigroup on $\mathcal{X}$. In this case $\rho(\mathcal{G}) \neq \emptyset$, hence the previous result gives a complete description of the spectrum and its finer subdivisions.

### 4.7. Spectral Theory for Complete Second Order Cauchy Problems.

We now apply our results to the reduction matrix
$$\mathcal{G} := \begin{pmatrix} 0 & Id \\ A & P \end{pmatrix} \tag{4.11}$$
associated to the complete second order Cauchy problem
$$\ddot{u}(t) = P\dot{u}(t) + Au(t). \tag{ACP_2}$$
We note that only in case $P = 0$ there is a satisfactory theory for (ACP$_2$), see, e.g., [3, Sect. 3.14]. In the complete case, i.e. if $P \neq 0$, there are many partial results and we refer to [7, Sect. VI.2] for a review of some of them.

Here we consider the following setting. For the definition of the extrapolated operator $A^Z$ on the extrapolation space $Z_{-1}$ see Proposition A.2.

- $U$, $X$, $Z$ are Banach spaces and $\mathcal{X} := Z \times X$,
- $A : D(A) \subseteq X \to X$ is an operator on $X$ satisfying Assumption A.1,

- $[D(A)] \hookrightarrow Z \hookrightarrow X$,
- $P = BC \in \mathcal{L}(Z, Z_{-1})$ where $B \in \mathcal{L}(U, Z_{-1})$ and $C \in \mathcal{L}(Z, U)$.

Under these hypotheses we consider the operator matrix $\mathcal{G} : D(\mathcal{G}) \subseteq \mathcal{X} \to \mathcal{X}$ in (4.11) equipped with the domain
$$D(\mathcal{G}) := \left\{ \begin{pmatrix} z \\ x \end{pmatrix} \in Z \times Z : A^Z z + BCx \in X \right\}.$$
As a first step towards the description of the spectrum of $\mathcal{G}$ we represent it as $\mathcal{G} = \mathcal{A}_{\mathcal{BC}}$. Define $\mathcal{Z} := Z \times Z \hookrightarrow \mathcal{X} = Z \times X \hookrightarrow \mathcal{Z}_{-1} := Z \times Z_{-1}$ and $\mathcal{U} := Z \times U$. Moreover, for some fixed $\mu_0 \in \rho(A)$ consider the operators
$$\mathcal{A}^{\mathcal{Z}} := \begin{pmatrix} 0 & Id \\ A^Z - \mu_0 & 0 \end{pmatrix} : \mathcal{Z} \subseteq \mathcal{Z}_{-1} \to \mathcal{Z}_{-1} \quad \text{and}$$
$$\mathcal{B} := \begin{pmatrix} 0 & 0 \\ Id & B \end{pmatrix} \in \mathcal{L}(\mathcal{U}, \mathcal{Z}_{-1}), \quad \mathcal{C} := \begin{pmatrix} \mu_0 & 0 \\ 0 & C \end{pmatrix} \in \mathcal{L}(\mathcal{Z}, \mathcal{U}).$$
Then $0 \in \rho(\mathcal{A}^{\mathcal{Z}})$ with resolvent

(4.12) $$R(0, \mathcal{A}^{\mathcal{Z}}) = \begin{pmatrix} 0 & R(\mu_0, A^Z) \\ -Id & 0 \end{pmatrix}.$$

In particular, $\rho(\mathcal{A}^{\mathcal{Z}}) \neq \emptyset$, hence the Assumptions 1.1 are satisfied. Let
$$\mathcal{G}^{\mathcal{Z}} := \mathcal{A}^{\mathcal{Z}} + \mathcal{BC} = \begin{pmatrix} 0 & Id \\ A^Z & BC \end{pmatrix} : \mathcal{Z} \subseteq \mathcal{Z}_{-1} \to \mathcal{Z}_{-1},$$
then $\mathcal{G}^{\mathcal{Z}} \begin{pmatrix} z \\ x \end{pmatrix} \in \mathcal{X}$ iff $\begin{pmatrix} z \\ x \end{pmatrix} \in D(\mathcal{G})$ which proves $\mathcal{G} = \mathcal{G}^{\mathcal{Z}}|_{\mathcal{X}} = \mathcal{A}_{\mathcal{BC}}$. Using this representation we obtain the following result, where for $\lambda \in \mathbb{C}$ we put
$$Q(\lambda) := \lambda^2 - \lambda \cdot BC - A^Z : Z \subseteq Z_{-1} \to Z_{-1}.$$

**Corollary 4.13.** *Let $\lambda \in \mathbb{C}$. Then*
$$\lambda \in \sigma_p(\mathcal{G}) \iff 0 \in \sigma_p(Q(\lambda)).$$
*Moreover, if there exists $\nu \in \mathbb{C}$ such that $0 \in \rho(Q(\nu))$, then*
$$\lambda \in \sigma(\mathcal{G}) \iff 0 \in \sigma(Q(\lambda)),$$
$$\lambda \in \sigma_*(\mathcal{G}) \iff 0 \in \sigma_*(Q(\lambda))$$
*for all $* \in \{a, r, c, ess\}$.*

*Proof.* We apply Corollary 2.7 for $\lambda_0 = 0 \in \rho(\mathcal{A}^{\mathcal{Z}})$ and $W = \mathcal{Z}$. To this end we first compute
$$\Delta_{\mathcal{Z}}(0, \lambda) = R(0, \mathcal{A}^{\mathcal{Z}}) \cdot (Id_{\mathcal{X}}, \mathcal{B}) \cdot \begin{pmatrix} -\lambda \cdot Id_{\mathcal{Z}} \\ \mathcal{C} \end{pmatrix}$$
$$= \begin{pmatrix} 0 & R(\mu_0, A^Z) \\ -Id & 0 \end{pmatrix} \cdot \begin{pmatrix} Id & 0 & 0 & 0 \\ 0 & Id & Id & B \end{pmatrix} \cdot \begin{pmatrix} -\lambda & 0 \\ 0 & -\lambda \\ \mu_0 & 0 \\ 0 & C \end{pmatrix}$$
$$= \begin{pmatrix} \mu_0 R(\mu_0, A^Z) & R(\mu_0, A^Z) \cdot (BC - \lambda) \\ \lambda & 0 \end{pmatrix} \in \mathcal{L}(\mathcal{Z}),$$
where we applied (4.12) and (2.10). Using Corollary 2.7.(a), Schur complements from Lemma A.9.(i) and the resolvent equation we conclude
$$\lambda \in \sigma_p(\mathcal{G}) \iff 1 \in \sigma_p(\Delta_{\mathcal{Z}}(0, \lambda))$$
$$\iff Id - \mu_0 R(\mu_0, A^Z) - \lambda R(\mu_0, A^Z) \cdot (BC - \lambda)$$
$$= R(\mu_0, A^Z) \cdot Q(\lambda) : Z \to Z \quad \text{is injective}$$
$$\iff 0 \in \sigma_p(Q(\lambda)).$$

Finally, if $0 \in \rho(Q(\nu))$ for some $\nu \in \mathbb{C}$, then by a similar reasoning we conclude $1 \in \rho(\Delta_{\mathcal{Z}}(0, \nu))$ and the remaining assertion follows from Corollary 2.7.(c). □

## 5. Conclusion

Our main result Theorem 2.3 characterizes spectral values of operators $A_{BC} := (A^Z + BC)|_X$ for triples $(A^Z, B, C) \in \mathcal{L}(Z, Z_{-1}) \times \mathcal{L}(U, Z_{-1}) \times \mathcal{L}(Z, U)$ which might perturb both the action and the boundary conditions of an operator $A = A^Z|_X$ on a Banach space $X$. Due to its generality this allows to study systematically and in a unified way spectral properties of various operators, thus furnishing an important tool for the analysis of the asymptotic behavior of the associated abstract Cauchy problem (ACP). For our results we impose only minimal assumptions, in particular we do not rely on admissibility conditions, denseness or closedness assumptions or a Hilbert space structure as used in [6, 16, 19]. Our approach is based on an extrapolated operator $A_{BC}^Z : Z \subset Z_{-1} \to Z_{-1}$ to which we associate an operator matrix having essentially the same spectral properties. The invertibility of this matrix is then investigated by using Schur complements which gives information on the spectral values of $A_{BC}^Z$. Finally, we use a result on the spectrum of the part of an operator to return to $A_{BC} = A_{BC}^Z|_X$.

## Appendix A.

In this appendix we introduce our notation, provide some results on extrapolated operators, consider spectral properties of parts of operators and study the invertibility of operator matrices by means of so-called Schur complements.

A.1. **Notation.** Besides the common notions for sets of numbers we use the abbreviations $\mathbb{N} := \{1, 2, 3, \ldots\}$, $\mathbb{N}_0 := \{0, 1, 2, \ldots\}$ and $\mathbb{R}_+ := [0, +\infty)$. By $(X, \|\bullet\|_X)^\sim$ we intend the completion of the normed space $X$. For normed spaces $X$ and $Y$ we denote by $\mathcal{L}(X, Y)$ the normed space of all bounded linear operators from $X$ to $Y$ and set $\mathcal{L}(X) := \mathcal{L}(X, X)$. Moreover, by $X \hookrightarrow Y$ we indicate a continuous embedding of $X$ in $Y$. For an operator $T : D(T) \subseteq X \to Y$ between Banach spaces $X$ and $Y$ we define $[D(T)] := (D(T), \|\bullet\|_T)$ for the *graph norm* given by $\|x\|_T := \|x\|_X + \|Tx\|_Y$. Then $[D(T)]$ is complete if and only if $T$ is closed. The transposed of a vector is denoted by $(\cdots)^\top$, while $\text{diag}(\cdots)$ indicates a diagonal matrix.

For notions related to the spectrum and resolvent of a linear operator, see Definition A.6.

A.2. **Abstract Extrapolation Spaces.** To apply our abstract results to a given operator $G$ we have to represent it as $G = (A^Z + BC)|_X$ for suitable operators $A^Z, B, C$ between spaces $Z, Z_{-1}, U$, cf. Diagram 1. Here the first step is usually to extend an operator $A$ with domain $D(A)$ on $X$ to $A^Z$ with bigger domain $Z \supset D(A)$ on a bigger space $Z_{-1} \supset X$.

In Section 3 we showed how this can be achieved in the context of our generic example. In this section we will introduce "abstract extrapolation spaces" to do so. This approach is more general as the construction of $Z_{-1}$ and $A^Z$ in Subsection 3.3 since it does not rely on a special form of the operators $B$ and $C$. However, it has the drawback that we will need some kind of denseness assumption on $A$, cf. Assumption A.1.

For a linear operator $A : D(A) \subseteq X \to X$ on a Banach space $X$ we define for $n \in \mathbb{N}$ the Banach spaces
$$X_n := \overline{D(A^n)}^{\|\bullet\|_X}$$
equipped with the norm induced by $X$. Moreover, we consider the operator
$$P_n := A^n|_{X_n} : D(P_n) \subseteq X_n \to X_n \quad \text{with domain} \quad D(P_n) := \{x \in D(A^n) : A^n x \in X_n\}.$$
To proceed we need to make the following

*Assumption* A.1. Suppose that $A : D(A) \subset X \to X$ satisfies
 (i) $\rho(A) \neq \emptyset$,
 (ii) there exists $n_0 \in \mathbb{N}$ such that $P_{n_0}$ is densely defined, i.e., for all $x \in X_{n_0}$ and $\varepsilon > 0$ there exists $z \in D(A^{n_0})$ such that $A^{n_0} z \in X_{n_0}$ and $\|x - z\|_X < \varepsilon$.

Under this assumption $A$ can be extended from $D(A)$ to a bigger domain $Z$ without changing its spectrum. More precisely, the following holds.

**Proposition A.2.** *Let Assumption A.1 be satisfied. Then for every Banach space $Z$ satisfying $[D(A)] \hookrightarrow Z \hookrightarrow X$ there exists a Banach space $Z_{-1}$ satisfying $X \hookrightarrow Z_{-1}$ and an "extrapolated" operator $A^Z : Z \subseteq Z_{-1} \to Z_{-1}$ such that $A^Z|_X = A$ and $\sigma(A) = \sigma(A^Z)$.*

*Proof.* Fix some $\eta_0 \in \rho(A)$ and define for $n \in \mathbb{N}$ the Banach spaces
$$X_{-n} := \bigl(X, \|\cdot\|_{-n}\bigr)^\sim \quad \text{where} \quad \|x\|_{-n} := \|R(\eta_0, A)^n x\|_X \quad \text{for } x \in X.$$
Let $x \in X$ and $\varepsilon > 0$. Then $R(\eta_0, A)^{n_0} x \in D(A^{n_0}) \subseteq X_{n_0}$, hence by Assumption A.1 there exists $z \in D(P_{n_0})$ such that
$$\|R(\eta_0, A)^{n_0} x - z\|_X = \|x - (\eta_0 - A)^{n_0} z\|_{-n_0} < \varepsilon.$$
Since $(\eta_0 - A)^{n_0} z \in X_{n_0}$ this shows that $X_{n_0} \subset (X, \|\cdot\|_{-n_0})$ is dense. Moreover, by definition $X$ is dense in $X_{-n_0}$, hence $X_{n_0}$ is dense in $X_{-n_0}$. Next from
$$X_{n_0} = \overline{D(A^{n_0})}^{\|\cdot\|_X} \subseteq \overline{D(A)}^{\|\cdot\|_{-n_0}} \subseteq X_{-n_0}$$
we conclude that also $D(A)$ is dense in $X_{-n_0}$. Moreover, a simple computation shows that
$$(\eta_0 - A) : D(A) \subseteq X_{-n_0} \to X_{-(n_0+1)}$$
is an isometry with dense range, hence admits a unique bounded extension
$$(\eta_0 - A_{-(n_0+1)}) \in \mathcal{L}(X_{-n_0}, X_{-(n_0+1)})$$
which is a surjective isometry, hence invertible. Now consider the Banach spaces
$$Z_{-1} := \bigl((\eta_0 - A_{-(n_0+1)})Z, \|\cdot\|_{Z_{-1}}\bigr) \quad \text{where} \quad \|z\|_{Z_{-1}} := \|(\eta_0 - A_{-(n_0+1)})^{-1} z\|_Z,$$
$$X_{-1} := \bigl((\eta_0 - A_{-(n_0+1)})X, \|\cdot\|_{-1}\bigr) \quad \text{where} \quad \|x\|_{-1} := \|(\eta_0 - A_{-(n_0+1)})^{-1} x\|_X,$$
and the operators
$$A^Z := A_{-(n_0+1)}|_{Z_{-1}} : Z \subseteq Z_{-1} \to Z_{-1},$$
$$A_{-1} := A_{-(n_0+1)}|_{X_{-1}} : X \subseteq X_{-1} \to X_{-1},$$
cf. Diagram 5. Since by construction $A_{-(n_0+1)}$ extends $A$ we have that $A^Z|_X = A$ and $A_{-1}|_Z = A^Z$.

Diagram 5: Extrapolation of a linear operator $A$.

Moreover, by definition $(\eta_0 - A_{-1}) : X \to X_{-1}$ is a surjective isometry and one easily verifies that
$$A = (\eta_0 - A_{-1})^{-1} \cdot A_{-1} \cdot (\eta_0 - A_{-1}).$$
Thus, $A$ and $A_{-1}$ are similar which implies $\rho(A_{-1}) = \rho(A) \neq \emptyset$. Finally, by Lemma A.7.(vii) applied to $F = X_{-1}$, $E = Z_{-1}$ and $T = A_{-1}$ we conclude that $\rho(A^Z) = \rho(A)$ as claimed. □

*Example* A.3. On $X := C[0,1]$ consider the operator $A := \frac{d}{ds}$ with (non-dense) domain $D(A) := \{f \in C^1[0,1] : \int_0^1 f(s)\,ds = 0\}$. Then $P_n$ is densely defined on $X_n$ if and only if $n \geq 2$, i.e., $A$ verifies Assumption A.1 for $n_0 = 2$ (but not for $n_0 = 1$).

Finally, we give a sufficient resolvent condition implying Assumption A.1 which is in particular satisfied for Hille–Yosida operators, cf. [15].

**Lemma A.4.** *If there exists a sequence $\lambda_n \in \rho(A)$ such that*
$$\lim_{n \to +\infty} R(\lambda_n, A)x = 0 \quad \text{for all } x \in X$$
*then Assumption A.1 is verified for $n_0 = 1$.*

*Proof.* Let $x \in D(A)$. Then the resolvent equation implies
$$x = x + \lim_{n \to +\infty} R(\lambda_n, A)Ax = \lim_{n \to +\infty} \lambda_n R(\lambda_n, A)x.$$
Since $\lambda_n R(\lambda_n, A)x \in D(A^2) \subseteq D(P_1)$ this shows that $D(P_1)$ is dense in $D(A)$. Moreover, $D(A)$ is dense in $X_1$ by definition and hence $D(P_1)$ is dense in $X_1$ as claimed. □

*Remark* A.5. (i) If $A$ is densely defined, the above construction can be simplified considerably. In particular, one can define $X_{-1} := (X, \|\bullet\|_{-1})^\sim$ for the norm $\|x\|_{-1} := \|R(\eta_0, A)x\|_X$ and then immediately extend $A$ continuously to $A_{-1} \in \mathcal{L}(X, X_{-1})$. For the details and more facts on inter- and extrapolation spaces as well as the associated abstract "Sobolev towers", see [7, Sect. II.5]
(ii) If a pair $(B, C) \in \mathcal{L}(U, X_{-1}) \times \mathcal{L}(Z, U)$ satisfies for some $\lambda_0 \in \rho(A)$ the condition
$$\text{rg}\big(R(\lambda_0, A_{-1})B\big) \subseteq Z = D(C), \tag{A.1}$$
then it is called *compatible* with respect to $A$. In this case, by the resolvent equation, (A.1) holds for all $\lambda_0 \in \rho(A)$. Moreover, (A.1) implies $\text{rg}(B) \subseteq Z_{-1}$ and by the closed graph theorem we obtain $B \in \mathcal{L}(U, Z_{-1})$. For results concerning the generator property of $A_{BC}$ for compatible pairs with respect to a generator $A$, see [1, 2].
(iii) Extrapolated spaces and operators are mainly used to give sense to calculations which are not defined a priori. The amazing fact is that in most cases one only needs to know their existence but not an explicit representation for them.

### A.3. Spectral Theory for Parts of Operators.
For convenience, we first recall the following notions from spectral theory.

**Definition A.6.** For a linear operator $T : D(T) \subset F \to F$ on a Banach space $F$ we define
$$\begin{aligned}
\sigma(T) &:= \big\{\lambda \in \mathbb{C} : \lambda - T \text{ is not invertible in } \mathcal{L}(F)\big\} & \text{(spectrum)}, \\
\rho(T) &:= \mathbb{C} \setminus \sigma(T) & \text{(resolvent set)}, \\
\sigma_p(T) &:= \{\lambda \in \mathbb{C} : \lambda - T \text{ is not injective}\} & \text{(point spectrum)}, \\
\sigma_a(T) &:= \left\{\lambda \in \mathbb{C} : \begin{array}{l}\lambda - T \text{ is not injective or}\\ \text{has non-closed range}\end{array}\right\} & \text{(approximative point spectrum)}, \\
\sigma_c(T) &:= \left\{\lambda \in \mathbb{C} : \begin{array}{l}\lambda - T \text{ is injective with}\\ \text{dense, non-closed range}\end{array}\right\} & \text{(continuous spectrum)}, \\
\sigma_r(T) &:= \left\{\lambda \in \mathbb{C} : \begin{array}{l}\lambda - T \text{ is injective with}\\ \text{non-dense range}\end{array}\right\} & \text{(residual spectrum)}, \\
\sigma_{\text{ess}}(T) &:= \left\{\lambda \in \mathbb{C} : \begin{array}{l}\dim(\ker(\lambda - T)) = +\infty \text{ or}\\ \text{codim}(\text{rg}(\lambda - T)) = +\infty\end{array}\right\} & \text{(essential spectrum)}.
\end{aligned}$$

Finally, for $\lambda \in \rho(A)$ we define the *resolvent operator* $R(\lambda, T) := (\lambda - T)^{-1} \in \mathcal{L}(F)$.

The next result generalizes [7, Props. IV.1.15 & IV.2.17] and connects some spectral properties of an operator $T$ on $F$ to those of its part $T|_E$ in a subspace $E$ of $F$.

**Lemma A.7.** *Let $T : D(T) \subset F \to F$ be a linear operator on a Banach space $F$, let $E$ be a Banach space satisfying $D(T) \subseteq E \hookrightarrow F$ and let*
$$T_1 := T|_E : D(T_1) \subseteq E \to E \quad \text{with domain} \quad D(T_1) := \{x \in D(T) : Tx \in E\}.$$
*Then the following holds.*

  (i) $\ker(T) = \ker(T_1)$; *in particular $T$ is injective $\iff T_1$ is injective.*
 (ii) $\mathrm{rg}(T_1) = \mathrm{rg}(T) \cap E$; *in particular $T$ is surjective $\implies T_1$ is surjective.*
(iii) $\mathrm{rg}(T)$ *is closed in $F$* $\implies \mathrm{rg}(T_1)$ *is closed in $E$.*
 (iv) $\mathrm{codim}(\mathrm{rg}(T)) < +\infty \implies \mathrm{codim}(\mathrm{rg}(T_1)) < +\infty$.
  (v) *If $E + \mathrm{rg}(T)$ is dense in $F$, then $\mathrm{rg}(T_1)$ is dense in $E$ $\implies \mathrm{rg}(T)$ is dense in $F$.*
 (vi) *If $T$ is closed, then $T_1$ is closed, and*
$$\rho(T) \subseteq \rho(T_1) \quad \text{and} \quad R(\lambda, T_1) = R(\lambda, T)|_E \quad \text{for all } \lambda \in \rho(T).$$
(vii) *If $\rho(T) \neq \emptyset$ and $F_1^T \hookrightarrow E$, then in (ii)–(v) always equivalence holds. In particular, in this case $\sigma(T) = \sigma(T_1)$.*

*Proof.* While one inclusion in both cases (i) and (ii) is clear, the respective other inclusion follows by the definition of $D(T_1)$ using the fact that $D(T) \subseteq E$.
To show (iii) take $y_n \in \mathrm{rg}(T_1)$ such that $y_n \to y \in E$ as $n \to +\infty$. Since $E \hookrightarrow F$ and $\mathrm{rg}(T)$ is closed in $F$, this implies $y \in \mathrm{rg}(T) \cap E = \mathrm{rg}(T_1)$, i.e., $\mathrm{rg}(T_1)$ is closed in $E$.
For (iv) assume that $\mathrm{codim}(\mathrm{rg}(T_1)) = +\infty$. Then there exists an infinite, linearly independent subset $S \subset E \setminus \mathrm{rg}(T_1)$. Since by (ii), $\mathrm{rg}(T_1) = \mathrm{rg}(T) \cap E$ we conclude $S \subset F \setminus \mathrm{rg}(T)$, i.e., $\mathrm{codim}(\mathrm{rg}(T)) = +\infty$.
To show (v) we assume that $\mathrm{rg}(T)$ is not dense in $F$. Then there exists $0 \neq \psi \in F'$ such that $\psi|_{\mathrm{rg}(T)} = 0$. Let $\varphi := \psi|_E \in E'$. If $\varphi = 0$, then $\psi|_{E+\mathrm{rg}(T)} = 0$ and by the denseness assumption it follows that $\psi = 0$ contradicting the choice of $\psi$. Hence, $\varphi \neq 0$ and $\varphi|_{\mathrm{rg}(T_1)} = 0$ which implies that $\mathrm{rg}(T_1)$ is not dense in $E$.
For (vi) take $x_n \in D(T_1)$ such that $x_n \to x \in E$ and $T_1 x_n \to y \in E$ as $n \to +\infty$. Since $E \hookrightarrow F$ this implies $x_n \to x$ in $F$ and $Tx_n \to y$ in $F$ as $n \to +\infty$. By the closedness of $T$ this gives $x \in D(T)$ and $Tx = y$. From $y \in E$ it follows that $x \in D(T_1)$ and $T_1 x = y$, i.e., $T_1$ is closed. Now take $\lambda \in \rho(T)$. Then $R := R(\lambda, T)|_E$ is a closed algebraic inverse of $\lambda - T_1$ defined on all of $E$ and having range in $E$. By the closed graph theorem this implies $R \in \mathcal{L}(E)$, i.e., $\lambda \in \rho(T_1)$ and $R = R(\lambda, T_1)$. This shows (i)–(vi).
To verify (vii) we first define
$$T_2 := T_1|_{F_1^T} : D(T_2) \subseteq F_1^T \to F_1^T \quad \text{with domain} \quad D(T_2) := \{x \in D(T_1) : T_1 x \in F_1^T\}.$$
Then the pair $T_2, T_1$ satisfies the assumptions made for $T_1, T$, hence we can repeat the reasoning in (ii)-(v) with $T_1, T$ replaced by $T_2, T_1$, respectively. For (v) note that for $\mu \in \rho(T) \subseteq \rho(T_1)$ we always have $E = \mathrm{rg}(\mu - T_1) \subseteq F_1^T + \mathrm{rg}\, T_1$, hence the denseness assumption is automatically satisfied. Moreover, for such $\mu$ the operator $\mu - T_1 \in \mathcal{L}(F_1^T, F)$ is an isomorphism which induces a similarity transformation between $T_2$ and $T$. This implies that $T_2$ is surjective/has closed range/has range with finite co-dimension/has dense range, respectively, if and only if $T$ has. Summing up, this shows equivalence in (ii)–(v) if $\rho(T) \neq \emptyset$. $\square$

*Remark* A.8. Without the denseness assumption on $E + \mathrm{rg}(T)$ the assertion in Lemma A.7.(v) does not hold. To see this take an operator $S : D(S) \subset E \to E$ with dense range. Then for a Banach space $G \neq \{0\}$ define $F := E \oplus G$ and the operator $T : D(T) \subseteq F \to F$ by $Tx := Sx$ for $x \in D(T) := D(S)$. Then $T_1 := T|_E = S$ has dense range in $E$ while $\mathrm{rg}(T) = \mathrm{rg}(S) \subseteq E$ is not dense in $F$. Clearly, in this case $E + \mathrm{rg}(T) = E$ is not dense in $F$. Note that in this example $T$ is closed on $F$ if $S$ is closed on $E$.

The following is the main result of this section.

**Corollary A.9.** *In the situation of Lemma A.7 the following relations hold.*

 (i) $\sigma_p(T_1) = \sigma_p(T)$.
(ii) $\sigma(T_1) \subseteq \sigma(T)$.

(iii) $\sigma_a(T_1) \subseteq \sigma_a(T)$.
(iv) $\sigma_c(T_1) \subseteq \sigma_c(T)$ if $E + \mathrm{rg}(T)$ is dense in $F$.
(v) $\sigma_r(T_1) \supseteq \sigma_r(T)$ if $E + \mathrm{rg}(T)$ is dense in $F$.
(vi) $\sigma_{ess}(T_1) \subseteq \sigma_{ess}(T)$.
(vii) If $\rho(T) \neq \emptyset$ and $F_1^T \hookrightarrow E$, then in (ii)–(vi) always equality holds.

*Proof.* All assertions follow easily from Definition A.6 and the previous lemma applied to $\lambda - T$ for $\lambda \in \mathbb{C}$ instead of $T$. For (iv) & (v) note that $E + \mathrm{rg}(\lambda - T)$ is independent of $\lambda \in \mathbb{C}$. □

**A.4. Schur Complements for Operator Matrices.** In this section we give conditions characterizing various spectral properties of an operator matrix. This yields to the notion of "Schur complement" which in a certain sense generalizes the concept of determinant of scalar matrices to matrices with non-commuting entries.

**Lemma A.10.** *For Banach spaces $E, F, G, H$ and linear operators $P \in \mathcal{L}(E, G)$, $Q \in \mathcal{L}(F, G)$, $R \in \mathcal{L}(E, H)$, $S \in \mathcal{L}(F, H)$ define the operator matrix*

$$\mathcal{T} := \begin{pmatrix} P & Q \\ R & S \end{pmatrix} \in \mathcal{L}(E \times F, G \times H).$$

*Then the following holds.*

(i) *If $S \in \mathcal{L}(F, H)$ is invertible then for $\Delta_1 := P - QS^{-1}R \in \mathcal{L}(E, G)$ we have*

$$(A.2) \qquad \mathcal{T} = \begin{pmatrix} Id_G & QS^{-1} \\ 0 & Id_H \end{pmatrix} \cdot \begin{pmatrix} \Delta_1 & 0 \\ 0 & S \end{pmatrix} \cdot \begin{pmatrix} Id_E & 0 \\ S^{-1}R & Id_F \end{pmatrix}.$$

*Hence,*

$\mathcal{T} \in \mathcal{L}(E \times F, G \times H)$ *is injective/surjective/has closed range/has dense range, resp.*
$\iff \Delta_1 \in \mathcal{L}(E, G)$ *is injective/surjective/has closed range/has dense range, resp.*

*In particular, $\mathcal{T}$ is invertible iff $\Delta_1$ is invertible and in this case*

$$\mathcal{T}^{-1} = \begin{pmatrix} \Delta_1^{-1} & -\Delta_1^{-1} \cdot QS^{-1} \\ -S^{-1}R \cdot \Delta_1^{-1} & S^{-1} + S^{-1}R \cdot \Delta_1^{-1} \cdot QS^{-1} \end{pmatrix} \in \mathcal{L}(G \times H, E \times F).$$

*Moreover, $\dim(\ker(\mathcal{T})) = \dim(\ker(\Delta_1))$ and $\mathrm{codim}(\mathrm{rg}(\mathcal{T})) = \mathrm{codim}(\mathrm{rg}(\Delta_1))$.*

(ii) *If $P \in \mathcal{L}(E, G)$ is invertible then for $\Delta_2 := S - RP^{-1}Q \in \mathcal{L}(F, H)$ we have*

$$(A.3) \qquad \mathcal{T} = \begin{pmatrix} Id_G & 0 \\ RP^{-1} & Id_H \end{pmatrix} \cdot \begin{pmatrix} P & 0 \\ 0 & \Delta_2 \end{pmatrix} \cdot \begin{pmatrix} Id_E & P^{-1}Q \\ 0 & Id_F \end{pmatrix}.$$

*Hence,*

$\mathcal{T} \in \mathcal{L}(E \times F, G \times H)$ *is injective/surjective/has closed range/has dense range, resp.*
$\iff \Delta_2 \in \mathcal{L}(F, H)$ *is injective/surjective/has closed range/has dense range, resp.*

*In particular, $\mathcal{T}$ is invertible iff $\Delta_2$ is invertible and in this case*

$$\mathcal{T}^{-1} = \begin{pmatrix} P^{-1} + P^{-1}Q \cdot \Delta_2^{-1} \cdot RP^{-1} & -P^{-1}Q \cdot \Delta_2^{-1} \\ -\Delta_2^{-1} \cdot RP^{-1} & \Delta_2^{-1} \end{pmatrix} \in \mathcal{L}(G \times H, E \times F).$$

*Moreover, $\dim(\ker(\mathcal{T})) = \dim(\ker(\Delta_2))$ and $\mathrm{codim}(\mathrm{rg}(\mathcal{T})) = \mathrm{codim}(\mathrm{rg}(\Delta_2))$.*

*If $P$ and $S$ are both invertible, then the following holds.*

(iii) $\ker(\Delta_1) = P^{-1}Q \ker(\Delta_2)$ *and* $\ker(\Delta_2) = S^{-1}R \ker(\Delta_1)$.
(iv) $\Delta_1$ *is injective/surjective/has closed range/has dense range* $\iff$
    $\Delta_2$ *is injective/surjective/has closed range/has dense range, respectively.*
(v) $\dim(\ker(\Delta_1)) = \dim(\ker(\Delta_2))$ *and* $\mathrm{codim}(\mathrm{rg}(\Delta_1)) = \mathrm{codim}(\mathrm{rg}(\Delta_2))$.
(vi) $\Delta_1$ *is invertible if and only if $\Delta_2$ is invertible and in this case*

$$\Delta_1^{-1} = P^{-1} + P^{-1}Q \cdot \Delta_2^{-1} \cdot RP^{-1} \in \mathcal{L}(G, E),$$
$$\Delta_2^{-1} = S^{-1} + S^{-1}R \cdot \Delta_1^{-1} \cdot QS^{-1} \in \mathcal{L}(H, F).$$

*Proof.* (i)–(v) are simple consequences of the factorizations of $\mathcal{T}$ given in (A.2) and (A.3) using the fact that the upper/lower triangular matrices involved are all isomorphisms. The boundedness of the inverses of $\mathcal{T}$, $\Delta_1$ and $\Delta_2$ follows from the closed graph theorem. (vi) follows from (i) and (ii) by comparing the diagonal entries of the representations of $\mathcal{T}^{-1}$. □

*Remark* A.11. The operators $\Delta_1 = P - QS^{-1}R : E \to G$ and $\Delta_2 = S - RP^{-1}Q : F \to H$ appearing above are frequently called *Schur complements* of the matrix $\mathcal{T}$, cf. [12], [18, Defs. 1.6.1 & 2.2.12].

The previous result has the following useful application.

**Corollary A.12.** *Let $E$, $F$ be Banach spaces and $Q \in \mathcal{L}(F, E)$, $R \in \mathcal{L}(E, F)$. Then*
$$1 \in \sigma(QR) \iff 1 \in \sigma(RQ), \qquad 1 \in \sigma_*(QR) \iff 1 \in \sigma_*(RQ)$$
*for all $* \in \{p, a, r, c, ess\}$. Moreover, $\ker(Id_E - QR) = Q\ker(Id_F - RQ)$ and $\ker(Id_F - RQ) = R\ker(Id_E - QR)$. Finally, if $1 \in \rho(RQ)$ or, equivalently, $1 \in \rho(QR)$, then*

$$(A.4) \quad \begin{aligned} (Id_E - QR)^{-1} &= Id_E + Q(Id_F - RQ)^{-1}R, \\ (Id_F - RQ)^{-1} &= Id_F + R(Id_E - QR)^{-1}Q. \end{aligned}$$

*Proof.* In the situation of Lemma A.10 choose $G = E$, $H = F$, $P = Id_E$ and $S = Id_F$. Then $\Delta_1 = Id_E - QR$ and $\Delta_2 = Id_F - RQ$. Hence, all assertions concerning the spectra follow easily from the characterizations of the corresponding spectral properties (cf. Definition A.6) of $\mathcal{T}$, $\Delta_1$ and $\Delta_2$ in Lemma A.10.(iii)–(v). Finally, (A.4) follows from Lemma A.10.(vi). □

Martin Adler, Arbeitsbereich Funktionalanalysis, Mathematisches Institut, Auf der Morgenstelle 10, D-72076 Tübingen
*E-mail address*: maad@fa.uni-tuebingen.de

Klaus-Jochen Engel[1], Università degli Studi dell'Aquila, Dipartimento di Ingegneria e Scienze dell'Informazione e Matematica (DISIM), Via Vetoio, I-67100 L'Aquila (AQ)
*E-mail address*: klaus.engel@univaq.it

---

[1]Corresponding author.